\documentclass{article}

\usepackage[margin=1in]{geometry}
\usepackage[english]{babel}
\usepackage[utf8x]{inputenc}
\usepackage{comment}
\usepackage[T1]{fontenc}
\usepackage[pdftex]{graphicx}
\usepackage{subcaption}
\usepackage{amsthm,amsmath,amssymb}
\usepackage[colorlinks,citecolor=blue]{hyperref}
\usepackage[square,numbers]{natbib}
\usepackage{booktabs,float,multirow}
\usepackage{enumitem}
\usepackage[normalem]{ulem}
\usepackage{tablefootnote}
\usepackage[ruled,vlined,linesnumbered]{algorithm2e}
\DontPrintSemicolon
\usepackage{standalone}
\usepackage[toc,page,header]{appendix}
\usepackage{xpatch}
\xpretocmd{\appendixpagename}{\LARGE}{}{}

\usepackage{braket}
\usepackage{xspace}
\usepackage{tikz}
\usetikzlibrary{positioning,
                patterns,
                matrix,
                backgrounds
                }

\DeclareMathOperator{\GL}{GL} 

\DeclareMathOperator*{\argmin}{arg\,min}
\DeclareMathOperator*{\argmax}{arg\,max}

\newcommand{\bbC}{\mathbb{C}}	
	\newcommand{\bbF}{\mathbb{F}}

	\newcommand{\bbR}{\mathbb{R}}
\newcommand{\bbS}{\mathbb{S}}

\newcommand{\mcI}{\mathcal{I}}	\newcommand{\mcJ}{\mathcal{J}}
	
	\newcommand{\mcN}{\mathcal{N}}

\newcommand{\bfA}{\mathbf{A}}	\newcommand{\bfB}{\mathbf{B}}
	
\newcommand{\bfE}{\mathbf{E}}	\newcommand{\bfF}{\mathbf{F}}
	\newcommand{\bfH}{\mathbf{H}}
\newcommand{\bfI}{\mathbf{I}}	
\newcommand{\bfK}{\mathbf{K}}	\newcommand{\bfL}{\mathbf{L}}
\newcommand{\bfM}{\mathbf{M}}	\newcommand{\bfN}{\mathbf{N}}
\newcommand{\bfO}{\mathbf{O}}	
\newcommand{\bfQ}{\mathbf{Q}}	\newcommand{\bfR}{\mathbf{R}}
	\newcommand{\bfT}{\mathbf{T}}
\newcommand{\bfU}{\mathbf{U}}	\newcommand{\bfV}{\mathbf{V}}
\newcommand{\bfW}{\mathbf{W}}	\newcommand{\bfX}{\mathbf{X}}
\newcommand{\bfY}{\mathbf{Y}}	\newcommand{\bfZ}{\mathbf{Z}}

	\newcommand{\bfn}{\mathbf{n}}
	
\newcommand{\bfq}{\mathbf{q}}	
	
\newcommand{\bfu}{\mathbf{u}}	
	\newcommand{\bfx}{\mathbf{x}}

\newcommand{\balpha}{\boldsymbol{\alpha}}
\newcommand{\bbeta}{\boldsymbol{\beta}}
\newcommand{\bgamma}{\boldsymbol{\gamma}}

\newcommand{\btau}{\boldsymbol{\tau}}

\DeclareMathAlphabet\mathbfcal{OMS}{cmsy}{b}{n}

\newcommand{\mbcA}{\mathbfcal{A}}

\newcommand{\mbcQ}{\mathbfcal{Q}}	
	\newcommand{\mbcT}{\mathbfcal{T}}

\newcommand{\mbcY}{\mathbfcal{Y}}

	\newcommand{\mfR}{\mathfrak{R}}

\allowdisplaybreaks

\newtheorem{theorem}{Theorem}

\newtheorem{lemma}{Lemma}

\theoremstyle{definition}
\newtheorem{definition}{Definition}

\theoremstyle{remark}
\newtheorem{remark}{Remark}

\newcommand{\blostr}{\texttt{BLOSTR}\xspace}

\newcommand{\ketbra}[1]{|#1\rangle\langle#1|}
\newcommand{\brho}{\boldsymbol{\rho}}
\renewcommand{\i}{\mathrm{i}}
\renewcommand{\triangleq}{\equiv}

\title{A Provably Efficient Method for Tensor Ring Decomposition and Its Applications}

\author{
    Han Chen\thanks{Department of Statistical Science, Duke University. Supported by NIH R01HL169347.}
    \and
    Sitan Chen\thanks{Department of Computer Science, Harvard University.
    Supported by NSF CCF-2430375.}
    \and
    Anru R. Zhang\thanks{Department of Biostatistics \& Bioinformatics and
    Department of Computer Science, Duke University. Supported by NSF CAREER-2203741 and NIH R01HL169347.}}
\date{}

\begin{document}

\maketitle

\begin{abstract}
We present the first deterministic, finite-step algorithm for exact tensor ring (TR) decomposition. Our method leverages blockwise simultaneous diagonalization to recover TR cores from a limited number of tensor observations, under a dimension condition requiring each mode size to be at least quadratic in the TR rank and under a genericity assumption on the cores, thereby providing both algebraic insight and practical efficiency. We extend the approach to the symmetric TR setting, where parameter complexity is significantly reduced and applications arise naturally in physics-based modeling and exchangeable data analysis. To handle noisy observations, we develop a robust recovery scheme that couples our initialization with alternating least squares, achieving faster convergence and improved accuracy compared to classical methods. As applications, we obtain new algorithms for questions in other domains where tensor ring decomposition is a key primitive, namely matrix product state tomography in quantum information and provable learning of pushforward distributions in the foundations of machine learning. These contributions advance the algorithmic foundations of TR decomposition and open new opportunities for scalable tensor network computation.
\end{abstract}

\section{Introduction}\label{sec:intro}

Tensors, or multidimensional arrays, are ubiquitous in the natural sciences and engineering. Tensor decomposition techniques~\cite{acar2008unsupervised,de2009survey,fu2020computing,kolda2009tensor,papalexakis2016tensors,sidiropoulos2017tensor} are powerful tools for handling complex data, with applications in machine learning~\cite{ji2019survey,kossaifi2020tensor,panagakis2021tensor,rabanser2017introduction}, signal processing~\cite{cichocki2015tensor}, and neuroscience~\cite{cong2015tensor,liu2022characterizing,morup2006parallel,sedighin2024tensor}. Despite their wide use for capturing multilinear structure, classical models such as CP~\cite{carroll1970analysis,harshman1970parafac} and Tucker~\cite{tucker1966some} can face scalability issues for higher-order tensors, due to rapidly growing parameter counts and the difficulty of selecting appropriate ranks~\cite{kolda2009tensor}.

Tensor networks (TNs) address these bottlenecks by representing a high-order tensor as a network of interconnected low-order tensors, where contractions over shared indices encode structural dependencies~\cite{bridgeman2017hand,orus2014practical}. Such representations can yield substantially lower effective ranks than classical notions~\cite{bernardi2023dimension,landsberg2012geometry,ye2018tensor}. In particular, the Tensor Train (TT) decomposition~\cite{oseledets2011tensor} expresses an order-$d$ tensor as a sequence of $d$ three-way \emph{cores} connected linearly, reducing storage from exponential to linear in $d$. The Tensor Ring (TR) decomposition~\cite{mickelin2020algorithms,zhao2016tensor} closes the chain into a ring by connecting the first and last cores via a trace, removing boundary-rank constraints and inducing a circularly invariant representation of modes. 

Formally, for any field $\bbF$, an order-$d$ tensor $\mbcT\in\bbF^{n_1\times\cdots\times n_d}$ admits a TR decomposition if there exist cores $\mbcQ_k\in\bbF^{n_k\times r\times r}$ such that
\begin{equation}\label{eq:order4TR}
  T(\alpha_1,\ldots,\alpha_d)
  = \mathrm{tr}\!\left\{\bfQ_1^{(\alpha_1)}\bfQ_2^{(\alpha_2)}\cdots\bfQ_d^{(\alpha_d)}\right\},
  \quad \text{for all}\ \,(\alpha_1,\ldots,\alpha_d)\in [n_1]\times\cdots\times[n_d],
\end{equation}
where $\bfQ_k^{(\alpha_k)}\in\bbF^{r\times r}$ denotes the mode-1 slice $(\mbcQ_k)_{\alpha_k,:,:}$ of the $k$-th three-way tensor, $\mbcQ_k$. We call $\mbcQ_k$ the $k$-th \emph{core}, and $r$ the \emph{TR rank} (also called the \emph{bond dimension} in the physics literature, e.g. \cite{bridgeman2017hand,orus2014practical}). We denote the TR decomposition by $\mbcT=\mfR(\mbcQ_1,\ldots,\mbcQ_d)$.

A number of algorithmic advances make the TR representation practical. The original work~\cite{zhao2016tensor} introduced SVD- and iteration-based routines that assume full tensor access. Alternating least-squares (TR-ALS) methods update each core via least squares and can handle noise or partial observations, though they are nonconvex and sensitive to initialization~\cite{chen2020tensorring,khoo2021efficient,wang2017efficient,wu2024nonlinear,yuan2019tensor}.

Beyond iterative schemes, several tensor models admit exact, finite-step constructions under suitable assumptions. For Tucker decomposition, truncating the higher-order SVD (HOSVD) yields an exact orthogonal Tucker representation in a fixed number of SVDs~\cite{kolda2009tensor}. For tensor train decomposition, TT-SVD provides a finite-step factorization~\cite{oseledets2011tensor}. For orthogonally decomposable CP tensors, spectral / simultaneous-diagonalization methods recover the factors exactly~\cite{kolda2001orthogonal}. Compared to Tucker, orthogonally decomposable CP, and TT decompositions, tensor rings possess a more intricate cyclic structure, making them substantially more challenging to analyze. To the best of our knowledge, no finite-step procedure for exact tensor ring decomposition has previously been established in the literature.

This paper proposes \blostr, a \emph{\underline{Blo}ckwise \underline{S}imultaneous diagonalization method for \underline{T}ensor \underline{R}ing decomposition}, which recovers the TR cores in a \emph{finite} number of steps by leveraging simultaneous diagonalization ideas~\cite{leurgans1993decomposition} adapted to the cyclic TR structure. \blostr operates on sparse observations through carefully chosen contractions that yield jointly diagonalizable matrices. Under the condition that $n_k \geq r^2$ 
for all $k \in [d]$ and that the TR cores are generic (i.e., drawn from a measure absolutely continuous with respect to the Lebesgue measure), we show that \blostr achieves exact TR decomposition with a sample complexity of $O(r^2 \sum_j n_j)$, matching the number of free parameters in the TR decomposition and thus attaining optimal sample efficiency. This sampling optimality is particularly valuable for high-dimensional or high-order tensors. Furthermore, we extend the method to the symmetric TR setting, which arises in quantum simulation and related TN applications~\cite{orus2014practical,schollwock2011density,verstraete2008matrix}, and we develop a robust variant that enables stable, approximate recovery under noise and missing data. 

The main contributions of this paper are as follows:
(1) We introduce a finite-step algorithm for exact TR decomposition under identifiable conditions, leveraging blockwise simultaneous diagonalization tailored to tensor rings.
(2) We extend this framework to symmetric TRs, yielding implications for quantum and tensor-network modeling.
(3) We develop a robust algorithm for TR decomposition in the presence of perturbations and demonstrate strong empirical performance from sparse and noisy observations, showing clear improvements over TR-ALS baselines~\cite{khoo2021efficient,wang2017efficient,yuan2019tensor}.
(4) We showcase applications of \blostr, including matrix product state tomography in quantum information and pushforward learning in high-dimensional statistics.

We emphasize that the condition \(n_k \ge r^2\) specifies the scope of the present analysis. This is the regime in which our blockwise simultaneous diagonalization argument yields exact recovery in a finite number of steps. The complementary regime \(n_k < r^2\) is substantially more delicate and is not addressed in this paper. 
The assumption \(n_k \ge r^2\) remains meaningful in applications. In the matrix product state setting, we discuss in Remark~\ref{remark:dimension-assumption-mps} why this condition is still relevant in practice. Likewise, in the polynomial transformation learning setting of Section~\ref{sec:main_moment}, the mode sizes correspond to output dimensions, which may naturally be large relative to the latent dimension \(r\).

The rest of this paper is organized as follows. Section~\ref{sec:preliminary} introduces notation, the TR definition, and key properties. Section~\ref{sec:procedure} presents \blostr and the exact recovery analysis. Section \ref{sec:symmetric} further expands the analysis to the symmetric case. Section~\ref{sec:robust} develops the robust algorithm and perturbation bounds. Section~\ref{sec:simulation} reports numerical experiments. Section~\ref{sec:mps} discusses connections to matrix product states and moment estimation. Section~\ref{sec:discussion} concludes.

\section{Notation and Preliminaries}\label{sec:preliminary}

We use boldface lowercase letters, boldface uppercase letters, and boldface calligraphic letters to denote vectors, matrices, and tensors, respectively. We use $[d]$ to denote the set $\{1,\dots,d\}$, with indices interpreted cyclically modulo~$d$; that is, we equip $[d]$ with the cyclic order $1 \prec 2 \prec \cdots \prec d \prec 1$, so that $d$ and~$1$ are consecutive, and addition or subtraction on $[d]$ is taken modulo~$d$. Here, ``$\prec$" denotes the precedence relationship. 
Let $\bfI_{d}$ be the identity matrix of size $d$, and $\mathbf{1}_d$ the all-ones vector in $\bbR^d$. For any field $\bbF$, let $\GL(n, \bbF)$ be the order-$n$ general linear group over $\bbF$. We use the notation $\bbF^{n^{\times d}}$ to denote the space of order-$d$ tensors with mode size $n$ over the field $\bbF$. In this paper, we focus on the case $\bbF=\bbC$, the field of complex numbers. More generally, the results remain valid over other algebraically closed fields or $\mathbb{R}$, subject to minor adjustments to the proposed algorithms. 

For a $d$-tuple $\balpha=(\alpha_1,\dots,\alpha_d)$ with $d\geq3$, let $\balpha_{\rm mid}$ be the $(d-2)$-tuple $(\alpha_2,\dots,\alpha_{d-1})$ that contains the middle $(d-2)$ elements of $\balpha$, and $\overleftarrow{\balpha}_{\rm mid}^k=(\alpha_{k+2},\dots,\alpha_d,\alpha_1,\dots,\alpha_{k-1})$ the middle $(d-2)$-tuple after circularly shifting $\balpha$ by $k$. We denote by $\Re(z)$ and $\Im(z)$ the real and imaginary parts of a complex number $z$, respectively. We denote the imaginary unit by $i=\sqrt{-1}$. For a matrix $\bfA\in\bbC^{m\times n}$, denote by $\bfA^\top$, $\bfA^*$ and $\bfA^\dagger$ the transpose, conjugate transpose and Moore-Penrose inverse of $\bfA$, respectively. Let $\text{tr}(\bfA)$ be the trace of $\bfA$. We use $\|\cdot\|_F$ to denote the Frobenius norm. 

The elements of a tensor of any order can be addressed by one of two ways; for example, the $(\alpha_1,\alpha_2,\alpha_3)$-th element of the tensor $\mbcT\in\bbC^{n_1\times n_2\times n_3}$ is addressed either $T(\alpha_1,\alpha_2,\alpha_3)$ or $T_{\alpha_1,\alpha_2,\alpha_3}$ for indices $\alpha_1\in[n_1]$, $\alpha_2\in[n_2]$, and $\alpha_3\in[n_3]$. We use $\bfA^{(\alpha)}\in\bbF^{r_1\times r_2}$ to denote the $\alpha$-th horizontal slice of the order-3 tensor $\mbcA\in\bbF^{n\times r_1\times r_2}$. For an order-$d$ tensor $\mbcT\in\bbF^{n_1\times\cdots\times n_d}$, let $\bfT_{[k]}, \bfT_{\left<k\right>}\in\bbF^{n_k\times \prod_{j\neq k}n_j}$ be its two matricizations:\footnote{ Note that $[\cdot]$ serves two distinct purposes in this paper: 
$[d] = \{1, \ldots, d\}$ denotes an index set, whereas $\mathbf{T}_{[k]}$ 
denotes the mode-$k$ matricization. The intended meaning is clear from 
context.} 
\begin{align}
    \label{eq:matricization}
    \bfT_{[k]}(\alpha_k,\overline{\alpha_{k+1}\cdots \alpha_d\alpha_1\cdots \alpha_{k-1}}) = \bfT_{\left<k\right>}(\alpha_k,\overline{\alpha_{k-1}\cdots \alpha_1\alpha_d\cdots \alpha_{k+1}})=T(\alpha_1,\cdots,\alpha_d),
\end{align}
where the vectorized index follows the column-major order: $\overline{\alpha_1 \alpha_2\cdots\alpha_d}\equiv \alpha_1+\sum_{j=2}^d(\alpha_j-1)\prod_{k=1}^{j-1}n_k$. The corresponding inverse operations transforming a matrix $\bfA\in\bbF^{n_k\times\prod_{j\neq k} n_j}$ to a tensor $\mbcA\in\bbF^{n_1\times n_2\times\cdots\times n_d}$ are 
\begin{align*}
    \left[\mathsf{Reshape}(\bfA, n_1,n_2,\dots, n_d,k)\right]_{\alpha_1,\alpha_2,\dots,\alpha_d} &= A(\alpha_k,\overline{\alpha_{k+1}\cdots \alpha_d\alpha_1\cdots \alpha_{k-1}}),\\
    [\mathsf{Reshape}'(\bfA,n_1, n_2,\dots, n_d,k)]_{\alpha_1,\alpha_2,\dots ,\alpha_d} &= A(\alpha_k,\overline{\alpha_{k-1}\cdots \alpha_1\alpha_d\cdots \alpha_{k+1}}),
\end{align*}
where $\alpha_k\in[n_k]$ for $k\in[d]$. For $d\geq3$ and $k\in[d]$, we define $\overleftarrow{\mbcT}^k\in\bbF^{n_{k+1}\times\cdots\times n_d\times n_1\times\cdots\times n_k}$ as circularly shifting the dimensions of $\mbcT$ by $k$. That is, $\overleftarrow{T}^k(\alpha_{k+1},\dots,\alpha_d,\alpha_1,\dots,\alpha_k)=T(\alpha_1,\dots,\alpha_k,\alpha_{k+1},\dots,\alpha_d)$.

For an index set $\Gamma \subseteq [n_d]$ of cardinality $r^2\leq n_d$, we use $\bfT( :, \balpha_{\rm mid}, \Gamma)=\bfT(:,\alpha_2,\dots,\alpha_{d-1},\Gamma)$ to denote the $n_1 \times r^2$ matrix obtained by fixing the mode-$2$ through mode-$(d-1)$ indices to $\alpha_2, \dots, \alpha_{d-1}$ and restricting the mode-$d$ index to $\Gamma$. In other words, this matrix consists of the mode-$d$ fibers $\{\bfT(:, \alpha_2, \dots, \alpha_{d-1}, \alpha_d)\}_{\alpha_d\in\Gamma}$ arranged as columns.

For a tensor $\mbcT \in \bbF^{n_1 \times \cdots \times n_d}$ and
$\bfU \in \bbF^{m \times n_k}$, the \emph{mode-$k$ product} $\mbcY = \mbcT \times_k \bfU \in \bbF^{n_1 \times \cdots \times n_{k-1} \times m \times n_{k+1} \times \cdots \times n_d}$
is defined entrywise as
\[
  Y(\alpha_1,\ldots,\alpha_{k-1}, j, \alpha_{k+1},\ldots,\alpha_d)
  \;=\;
  \sum_{\alpha_k=1}^{n_k} T(\alpha_1,\ldots,\alpha_k,\ldots,\alpha_d)\,
  U(j,\alpha_k).
\] 
Equivalently, in terms of the mode-$k$ matricization, $\bfY_{[k]} \;=\; \bfU \bfT_{[k]}.$

For two matrices $\bfA\in\bbF^{m\times n}$ and $\bfB\in\bbF^{p\times q}$, the \emph{Kronecker product} of $\bfA$ and $\bfB$ is 
\[
    \bfA \otimes \bfB =
    \begin{bmatrix}
    a_{11}\bfB & a_{12}\bfB & \cdots & a_{1n}\bfB \\
    a_{21}\bfB & a_{22}\bfB & \cdots & a_{2n}\bfB \\
    \vdots  & \vdots  & \ddots & \vdots  \\
    a_{m1}\bfB & a_{m2}\bfB & \cdots & a_{mn}\bfB
    \end{bmatrix}\in\bbF^{mp\times nq}.
\]
Equivalently, in terms of entries:
\[
    (\bfA \otimes \bfB)_{(j-1)p+r,\,(k-1)q+s} = a_{jk}b_{rs}, \text{for}\ j\in[m],\ k\in[n],\ r\in[p],\ \text{and}\ s\in[q].
\]

\begin{figure}[ht]
    \centering
    \begin{tikzpicture}[every node/.style={minimum size=1cm, anchor=center}]
    \definecolor{block1}{HTML}{FFFFFF}
    \definecolor{block2}{HTML}{F8F9FA}
    \definecolor{block3}{HTML}{DEE2E6}
    \definecolor{block4}{HTML}{ADB5BD}

    \matrix (X) [matrix of nodes,nodes={minimum size=.7cm},column sep=-\pgflinewidth, row sep=-\pgflinewidth, label=below:{$\mathbf{X}$}, right delimiter=\rbrack,left delimiter=\lbrack] {
        $X_{11}$ & $X_{12}$ & $X_{13}$ & $X_{14}$ & $X_{15}$ & $X_{16}$ \\
        $X_{21}$ & $X_{22}$ & $X_{23}$ & $X_{24}$ & $X_{25}$ & $X_{26}$ \\
        $X_{31}$ & $X_{32}$ & $X_{33}$ & $X_{34}$ & $X_{35}$ & $X_{36}$ \\
        $X_{41}$ & $X_{42}$ & $X_{43}$ & $X_{44}$ & $X_{45}$ & $X_{46}$ \\
        $X_{51}$ & $X_{52}$ & $X_{53}$ & $X_{54}$ & $X_{55}$ & $X_{56}$ \\
        $X_{61}$ & $X_{62}$ & $X_{63}$ & $X_{64}$ & $X_{65}$ & $X_{66}$ \\
    };

    \begin{pgfonlayer}{background}
        \fill[block1] (X-1-1.north west) rectangle (X-3-3.south east);
        \fill[block2] (X-1-4.north west) rectangle (X-3-6.south east);
        \fill[block3] (X-4-1.north west) rectangle (X-6-3.south east);
        \fill[block4] (X-4-4.north west) rectangle (X-6-6.south east);
    \end{pgfonlayer}

    \matrix (PX) [matrix of nodes,nodes={minimum size=.7cm},column sep=-\pgflinewidth, row sep=-\pgflinewidth, right=3cm of X, label=below:{$\Pi_{2,3}(\mathbf{X})$}, right delimiter=\rbrack,left delimiter=\lbrack] {
        $X_{11}$ & $X_{14}$ & $X_{12}$ & $X_{15}$ & $X_{13}$ & $X_{16}$ \\
        $X_{41}$ & $X_{44}$ & $X_{42}$ & $X_{45}$ & $X_{43}$ & $X_{46}$ \\
        $X_{21}$ & $X_{24}$ & $X_{22}$ & $X_{25}$ & $X_{23}$ & $X_{26}$ \\
        $X_{51}$ & $X_{54}$ & $X_{52}$ & $X_{55}$ & $X_{53}$ & $X_{56}$ \\
        $X_{31}$ & $X_{34}$ & $X_{32}$ & $X_{35}$ & $X_{33}$ & $X_{36}$ \\
        $X_{61}$ & $X_{64}$ & $X_{62}$ & $X_{65}$ & $X_{63}$ & $X_{66}$ \\
    };

    \begin{pgfonlayer}{background}
        \foreach \i in {1,...,6} {
            \foreach \j in {1,...,6} {
                \pgfmathtruncatemacro{\iOff}{mod(\i-1,2)}
                \pgfmathtruncatemacro{\jOff}{mod(\j-1,2)}
                \pgfmathtruncatemacro{\colorIndex}{\iOff*2 + \jOff + 1}
                \edef\colorname{block\colorIndex}
                \fill[\colorname] (PX-\i-\j.north west) rectangle (PX-\i-\j.south east);
            }
        }
    \end{pgfonlayer}

    \draw[->, thick, shorten <=25pt, shorten >=25pt] (X-3-6.south east) -- (PX-3-1.south west) node[midway, above] {$\Pi_{2,3}$};
\end{tikzpicture}
    \caption{Illustration of $\Pi_{2,3}$ operated on a 6-by-6 matrix $\bfX$. Elements in the same sub-matrices are labeled with same colors.}
    \label{fig:Pi}
\end{figure}
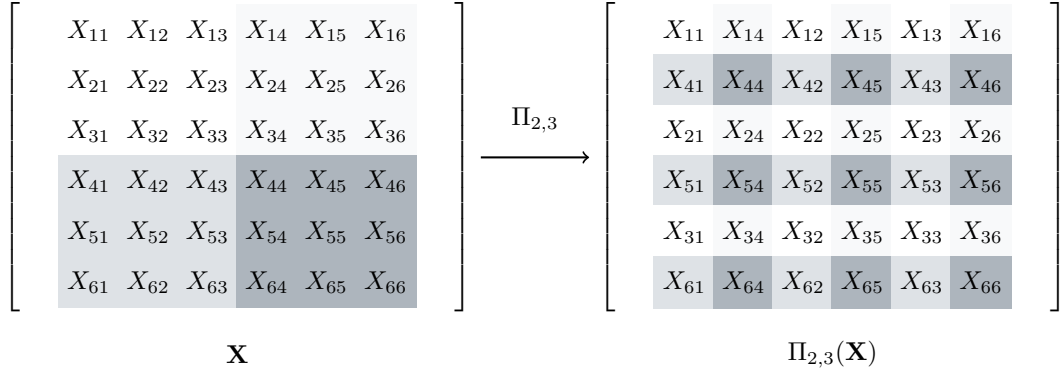

For integers $r_1, r_2\geq 1$, we define a row-column-wise permutation operation $\Pi_{r_1,r_2}(\cdot)$ on $r_1r_2$-by-$r_1r_2$ matrices: For a matrix $\bfX\in\bbF^{r_1r_2\times r_1r_2}$,
\begin{align}\label{eq:pi}
    \left(\Pi_{r_1,r_2}(\bfX)\right)_{(j_2-1)r_1+j_1,(k_2-1)r_1+k_1}=X_{(j_1-1)r_2+j_2,(k_1-1)r_2+k_2}\ \text{for}\ j_1,k_1\in[r_1],\ j_2,k_2\in[r_2].
\end{align}
For simplicity, we omit the subscripts when $r_1=r_2$. Figure~\ref{fig:Pi} gives an illustration where $r_1=2$ and $r_2=3$.

The individual blocks of $\mbcQ_1,\dots,\mbcQ_d$ of a TR decomposition are unidentifiable due to the \emph{gauge invariance}: For any $\bfL_k\in\GL(r, \mathbb{C})$,
	\begin{equation}\label{eq:order-3-invariance}
	    \text{tr}\left\{\bfQ_1^{(\alpha_1)}\bfQ_2^{(\alpha_2)}\cdots\bfQ_d^{(\alpha_d)}\right\} =  \text{tr}\left\{\left(\bfL_{1}^{-1} \bfQ_1^{(\alpha_1)} \bfL_2\right)\left(\bfL_2^{-1} \bfQ_2^{(\alpha_2)}\bfL_3\right)\cdots \left(\bfL_{d}^{-1} \bfQ_d^{(\alpha_d)}\bfL_1\right)\right\}.
	\end{equation}
We thus introduce the following equivalence class: For all $\bfL_k\in\text{GL}(r,\bbC)$, 
\begin{equation}\label{eq:gauge-invariance}
    \left\{\mbcQ_1, \mbcQ_2,\dots, \mbcQ_d\right\} \sim \left\{\mbcQ_1 \times_2 \bfL_{1}^{-1} \times_3 \bfL_2^*, \mbcQ_2 \times_2 \bfL_2^{-1} \times_3 \bfL_3^*, \dots, \mbcQ_d \times_2 \bfL_{d}^{-1} \times_3 \bfL_1^*\right\}.
\end{equation}
The primary goal of this work is to identify a series of representatives $(\hat{\mbcQ}_1,\hat{\mbcQ}_2,\dots,\hat{\mbcQ}_d)$ that are equivalent to the original TR cores.

\paragraph{Order-2 TR Decomposition} It is helpful to first explore the TR decomposition for order-2 tensors, i.e., matrices, before diving into more general scenarios. Suppose we observe a matrix $\bfT\in \mathbb{C}^{n_1\times n_2}$ that admits a TR decomposition: $$T(\alpha_1,\alpha_2) = \text{tr}\left\{\bfQ_1^{(\alpha_1)} \bfQ_2^{(\alpha_2)}\right\}=\bfQ_{1\left<1\right>}(\alpha_1,:)\left(\bfQ_{2\left[1\right]}(\alpha_2,:)\right)^*, \quad \alpha_1\in [n_1],\ \alpha_2\in[n_2].$$
Then $\bfT=\bfQ_{1\left<1\right>}\bfQ_{2\left[1\right]}^*$ is an $n_1$-by-$n_2$ matrix that forms a rank-$r^2$ decomposition. A representative pair ${\hat{\mbcQ}_1, \hat{\mbcQ}_2}$ can be obtained directly via SVD. The complete procedure for order-2 TR decomposition is provided in Algorithm~\ref{alg:order-2}. Moreover, compared to that of the higher-order tensors, the TR cores of the matrix $\bfT$ are only identifiable up to the rotation of a significantly larger freedom ($r^2\times r^2$), compared to the tensor ring order-3 or higher ($r\times r$, as described in \eqref{eq:order-3-invariance}): for any $\bfL\in \GL(r^2, \mathbb{C})$,
$$\{\mbcQ_1, \mbcQ_2\} \sim \left\{\mathsf{Reshape'}(\bfQ_{1\left<1\right>}\bfL, n_1, r, r, 1), \mathsf{Reshape}(\bfQ_{2[1]}\bfL^{*-1}, n_2, r, r, 1)\right\}.$$

\begin{algorithm}[ht]  
    \setcounter{algocf}{-1}
    \caption{Order-2 Tensor Ring Decomposition}
    \label{alg:order-2}
    \KwIn{Tensor $\mbcT\in\bbC^{n_1\times n_2}$,  TR-rank $r$}
    \KwOut{TR-cores $\hat{\mbcQ}_1\in\bbC^{n_1\times r\times r}$, $\hat{\mbcQ}_2\in\bbC^{n_2\times r\times r}$}
    
    Obtain the SVD, $\bfT=\bfU\boldsymbol{\Sigma}\bfV^*$, where $\bfU\in\bbC^{n_1\times n_1}$, $\bfV\in\bbC^{n_2\times n_2}$, and $\boldsymbol{\Sigma}\in\bbC^{n_1\times n_2}$\;
    
    \If{$n_1\geq r^2$}{
        Let $\hat{\bfQ}_{1\left<1\right>}=\bfU(:,1\!:\!r^2)$, and $\hat{\bfQ}_{2[1]}=\bfV\left(\boldsymbol{\Sigma}(1\!:\!r^2,:)\right)^*$\;
    }
    \Else
        {Set $\hat{\bfQ}_{1\left< 1\right>} = [\bfU\ \bfO]$, where $\bfO$ is the $n_1 \times (r^2-n_1)$ zero matrix\;
        
        Define $\widetilde{\boldsymbol{\Sigma}} \in \mathbb{C}^{r^2\times n_2}$ such that its first $n_1\times n_2$ block equals $\boldsymbol{\Sigma}$ and the rest are zeros\;
        
        Let $\hat{\bfQ}_{2[1]} = \bfV\,\widetilde{\boldsymbol{\Sigma}}^*$\;
    }
    Reshape $\hat{\bfQ}_1,\hat{\bfQ}_2$ into $\hat{\mbcQ}_1,\hat{\mbcQ}_2$
\end{algorithm}

\paragraph{Ranks of TR Decomposition}
    We assume that the TR rank is known \emph{a priori} in this paper. Note that when it is not, the ranks can be estimated using methods such as SVD-based algorithms (see, e.g., \cite{zhao2016tensor}).
    
\section{Finite-Step TR Decomposition via Blockwise Simultaneous Diagonalization}\label{sec:procedure}

This section presents \blostr, a finite-step, closed-form procedure for tensor ring (TR) decomposition. Recall that the TR decomposition model \eqref{eq:order4TR} can be grouped as
\begin{align}\label{eq:def-R}
    T(\alpha_1,\alpha_2,\dots,\alpha_d)=\text{tr}\left\{\bfQ_1^{(\alpha_1)}\bfR^{\balpha_{\rm mid}}\bfQ_d^{(\alpha_d)}\right\},\quad\text{where}\ \bfR^{\balpha_{\rm mid}}\equiv\bfQ_2^{(\alpha_2)}\cdots\bfQ_{d-1}^{(\alpha_{d-1})}.
\end{align}
The slice of $\mbcT$ at index $\balpha_{\rm mid}$ satisfies
\begin{align}\label{eq:trace}
    \bfT(:,\balpha_{\rm mid},:)=\bfQ_{1\left<1\right>}\left(\bfI_r\otimes\bfR^{\balpha_{\rm mid}}\right)\bfQ_{d\left[1\right]}^*.
\end{align}

The following lemma shows the connection between the eigendecomposition of the tensor and that of the TR cores, which plays a key role in the proposed procedure.

\begin{lemma}\label{lemma:eigenvalues}
    Suppose $\mbcT\in\bbC^{n_1\times n_2\times\cdots\times n_d}$ satisfies \eqref{eq:order4TR}, where $d\geq 3$, $r\geq 2$, and $n_k \geq r^2$ for all $k\in[d]$. Assume that the elements of $\mbcQ_k\in\mathbb{C}^{n_k\times r\times r}$ are randomly drawn from a measure $\mu$ that is absolutely continuous with respect to the standard Euclidean measure. Let $\balpha,\bbeta\in[n_2]\times\cdots\times[n_{d-1}]$ with $\balpha\neq\bbeta$ be two distinct middle index tuples in the sense of $\balpha_{\mathrm{mid}}$ defined in Section~\ref{sec:preliminary}. Assume $\bfR^{\balpha}\left(\bfR^{\bbeta}\right)^{-1}$ has an eigendecomposition $\bfU\boldsymbol{\Lambda}\bfU^{-1}$, where $\boldsymbol{\Lambda}\in\bbC^{r\times r}$ is diagonalizable and has distinct diagonal elements $\lambda_1,\dots,\lambda_{r}$. $\Gamma_{\balpha}, \Gamma_{\bbeta} \subseteq [n_d]$ are two different subsets of cardinality $r^2$. Then with probability one, the nonzero eigenvalues of $\bfT(:,\balpha,\Gamma_{\balpha})\bfT(:,\bbeta,\Gamma_{\bbeta})^\dagger$ are $\lambda_1,\dots,\lambda_{r}$, each with multiplicity $r$, and the corresponding eigenvectors are
    \begin{align*}
        \bfE=\bfQ_{1\left<1\right>}\left(\bfI_{r}\otimes\bfU\right)\Pi(\bfK)\quad\text{for some}\ 
        \bfK=\begin{bmatrix}
            \bfK_1 & & \\
            & \ddots & \\
            & & \bfK_r 
        \end{bmatrix}\in\bbC^{r^2\times r^2},
    \end{align*}
    where $\bfK_j\in\GL(r,\bbC)$, $j\in[r]$. Equivalently, restricting to the nonzero spectrum (modulo the nullspace)
    \begin{flalign*}
        \bfT(:,\balpha,\Gamma_{\balpha})\bfT(:,\bbeta,\Gamma_{\bbeta})^\dagger&=\bfE(\bfI_r\otimes\boldsymbol{\Lambda})\bfE^\dagger\\
        &=\left(
        \bfQ_{1\left<1\right>}\underbrace{\begin{bmatrix}
            \bfU & & \\
            & \ddots & \\
            & & \bfU
        \end{bmatrix}}_{r\ \text{copies}}\Pi\left(
        \begin{bmatrix}
            \bfK_1 & & \\
            & \ddots & \\
            & & \bfK_r 
        \end{bmatrix}
        \right) 
        \right)
        \underbrace{\begin{bmatrix}
            \boldsymbol{\Lambda} & & \\
            & \ddots & \\
            & & \boldsymbol{\Lambda}
        \end{bmatrix}}_{r\ \text{copies}} \\
        &\hspace{80pt} \cdot \left(
        \bfQ_{1\left<1\right>}\underbrace{\begin{bmatrix}
            \bfU & & \\
            & \ddots & \\
            & & \bfU
        \end{bmatrix}}_{r\ \text{copies}}\Pi\left(
        \begin{bmatrix}
            \bfK_1 & & \\
            & \ddots & \\
            & & \bfK_r 
        \end{bmatrix}
        \right)
        \right)^\dagger.
    \end{flalign*}
\end{lemma}

Lemma~\ref{lemma:eigenvalues} implies that we can recover the first core (up to gauge) by exploiting the $r\times r$ block structure in the eigendecomposition of $\bfT(:,\balpha,\Gamma_{\balpha})\bfT(:,\bbeta,\Gamma_{\bbeta})^\dagger$. 
We outline the steps below. 
\paragraph{Step 1 (two spectral probes)} Choose $\balpha$, $\balpha'$, $\bbeta$ and $\bbeta'$ such that
\begin{align}
    \label{eq:alpha-beta}
    \balpha\neq\bbeta,\ \balpha'\neq\bbeta',\ \text{and}\ \{\balpha,\bbeta\}\neq\{\balpha',\bbeta'\}.
\end{align}
Compute an eigenbasis $\bfE$ for the nonzero spectrum of $\bfT(:,\balpha,\Gamma_{\balpha})\bfT(:,\bbeta,\Gamma_{\bbeta})^\dagger$, and let $\bfE'$ be the eigenbasis of $\bfT(:,\balpha',\Gamma_{\balpha'})\bfT(:,\bbeta',\Gamma_{\bbeta'})^\dagger$. Let $\bfU\boldsymbol{\Lambda}\bfU^{-1}$ and $\bfV\boldsymbol{\Lambda}'\bfV^{-1}$ be the eigendecompositions of $\bfR^{\balpha}(\bfR^{\bbeta})^{-1}$ and $\bfR^{\balpha'}(\bfR^{\bbeta'})^{-1}$, respectively. We order eigenvectors by eigenvalues and fix a basis within each multiplicity-$r$ eigenspace; any remaining rotation is absorbed into the block-diagonal gauge matrices $\bfK, \bfK'$. After calculation, we can obtain
\begin{equation*}
\begin{split}
    \bfF\equiv\bfE^\dagger\bfE'= &  \Pi(\bfK^{-1})\left(\bfI_{r}\otimes\left(\bfU^{-1}\bfV\right)\right)\Pi(\bfK')\\
    = &         \Pi \left(\begin{bmatrix}
            \bfK_1^{-1} & & \\
            & \ddots & \\
            & & \bfK_r^{-1} 
        \end{bmatrix} \right)\underbrace{\begin{bmatrix}
            \bfU^{-1}\bfV & & \\
            & \ddots & \\
            & & \bfU^{-1}\bfV
        \end{bmatrix}}_{r\ \text{copies}}         \Pi\left(\begin{bmatrix}
            \bfK_1' & & \\
            & \ddots & \\
            & & \bfK_r' 
        \end{bmatrix}\right).
\end{split}
\end{equation*}
This means for any $j,k\in[r]$,
\begin{align}
    \label{eq:F}
    \bfF_{[0:(r-1)]\cdot r+j,[0:(r-1)]\cdot r+k}=\left(\bfU^{-1}\bfV\right)_{j,k}\cdot\bfK_{j}^{-1}\bfK'_{k}.
\end{align}
\paragraph{Step 2 (block identification / gauge fixing)} Equation~\eqref{eq:F} encodes the relation between the block-wise structure of $\bfF$, the product of two \emph{known} matrices, and the block multiplication in $\bfK$, an \emph{unknown} matrix that, by Lemma~\ref{lemma:eigenvalues}, determines $\mbcQ_1$ up to the gauge invariance in \eqref{eq:gauge-invariance}. We fix the gauge by setting $\hat{\bfK}_1=\bfI_r$ and normalizing $(\bfU^{-1}\bfV)_{1,\ell}=(\bfU^{-1}\bfV)_{\ell,1}=1$.
Then, from \eqref{eq:F}, 
\begin{align*}
    \hat{\bfK}_\ell=\bfF_{[0:(r-1)]\cdot r+1,[0:(r-1)]\cdot r+1}(\bfF_{[0:(r-1)]\cdot r+\ell,[0:(r-1)]\cdot r+1})^{-1}.
\end{align*}

\paragraph{Step 3 (recover $\mbcQ_1$)} As $\bfE=\bfQ_{1\langle1\rangle}(\bfI_r\otimes\bfU)\Pi(\bfK)$, we obtain
\begin{align}
    \label{eq:hatq-1}
    \hat{\bfQ}_{1\left<1\right>}=\bfE(\Pi(\hat{\bfK}))^{-1}.
\end{align}

\paragraph{Step 4 (circular permutation and remaining cores)}
We have the following property of \emph{circular dimensional permutation invariance}:
\begin{lemma}[{\cite[Theorem 2.1]{zhao2016tensor}}] 
    \label{lemma:permutation}
    Let $\mbcT\in\bbF^{n_1\times n_2\times\cdots\times n_d}$ be an order-$d$ tensor and its TR decomposition is given by $\mbcT=\mfR(\mbcQ_1,\mbcQ_2,\dots,\mbcQ_d)$, where $\bbF$ is any field. Recall that \(\overleftarrow{\mbcT}^k\) denotes the tensor obtained from \(\mbcT\) by a circular shift of its modes by $k$, as defined in Section~\ref{sec:preliminary}. Then we have $\overleftarrow{\mbcT}^k=\mfR(\mbcQ_{k+1},\dots,\mbcQ_d,\mbcQ_1,\dots,\mbcQ_k)$.
\end{lemma}

Choose $\bgamma=(\gamma_1,\dotsm\gamma_d)\in[n_1]\times\cdots\times[n_d]$. By Lemma~\ref{lemma:permutation}, for each $2\leq k\leq d$,
\begin{align*}
    \overleftarrow{\bfT}^{k-1}(:,\overleftarrow{\bgamma}_{\rm mid}^{k-1},\Gamma_{k-1})=\bfQ_{k\left<1\right>}\left(\bfI_r\otimes(\bfQ_{k+1}^{(\gamma_{k+1})}\cdots\bfQ_1^{(\gamma_1)}\cdots\bfQ_{k-2}^{(\gamma_{k-2})})\right)\bfQ_{(k-1)[1]}^*(:,\Gamma_{k-1}).
\end{align*}
Proceeding sequentially, for $2\leq k\leq d-1$ we define
\begin{align}
    \label{eq:Qk}
    \hat{\bfQ}_{k\left<1\right>}=\overleftarrow{\bfT}^{k-1}(:,\overleftarrow{\bgamma}_{\rm mid}^{k-1},\Gamma_{k-1})\left(\left(\bfI_r\otimes(\hat{\bfQ}_1^{(\gamma_1)}\cdots\hat{\bfQ}_{k-2}^{(\gamma_{k-2})})\right)\hat{\bfQ}_{(k-1)[1]}^*(:,\Gamma_{k-1})\right)^\dagger,
\end{align}
with the convention that the product in parentheses is $\bfI_r$ when $k=2$.
We can show that the recovered TR cores satisfy the following identities:
\begin{align}\label{eq:q-equivalent-main}
    \hat{\mbcQ}_k=\left\{\begin{aligned}
        &\mbcQ_1\times_2 \bfK_1^* \times_3 \{((\bfU^{-1}\bfV)_{1,1})^{-1}\bfU\bfW\}^*, &&\text{if}\ k=1;\\
        &\mbcQ_2\times_2 (\bfU^{-1}\bfV)_{1,1}(\bfU\bfW)^{-1} \times_3 \{\bfQ_3^{(\gamma_3)}\cdots\bfQ_d^{(\gamma_d)}(\bfK_1^*)^{-1}\}^*, &&\text{if}\ k=2;\\
        &\mbcQ_k\times_2\bfK_1^*(\bfQ_k^{(\gamma_k)}\cdots\bfQ_d^{(\gamma_d)})^{-1}\times_3 \{\bfQ_{k+1}^{(\gamma_{k+1})}\cdots\bfQ_d^{(\gamma_d)}(\bfK_1^*)^{-1}\}^*, &&\text{if}\ 3\leq k\leq d-1;\\
        &\mbcQ_d\times_2 \bfK_1^*(\bfQ_d^{(\gamma_d)})^{-1}\times_3 \{(\bfK_1^*)^{-1}\}^*, &&\text{if}\ k=d;
    \end{aligned}\right.
\end{align}
where $\bfW\equiv\text{diag}\left((\bfU^{-1}\bfV)_{[1:r],1}\right)\in\bbC^{r\times r}$ and thus invertible with probability one under the assumption that the elements of each $\mbcQ_k\in\mathbb{C}^{n_k\times r\times r}$ are randomly drawn from a measure that is absolutely continuous with respect to the standard Euclidean measure. Details of the analysis are provided in the Appendix \ref{sec:proof-theorem1}. Therefore, we obtain all TR cores and finish the algorithm.

We formalize the above discussions into the following theorem.

\begin{theorem}[\blostr achieves exact tensor ring  decomposition in finite steps]\label{thm:order-d}
    Suppose $\mbcT\in\bbC^{n_1\times n_2\times\cdots\times n_d}$ satisfies \eqref{eq:order4TR}, where $d\geq 3$, $r\geq 2$, and $n_k \geq r^2$ for all $k\in[d]$. Assume that the elements of $\mbcQ_k\in\mathbb{C}^{n_k\times r\times r}$ are randomly drawn from a measure $\mu$ that is absolutely continuous with respect to the standard Euclidean measure. Then with probability one, $\mbcQ_1, \mbcQ_2,\dots, \mbcQ_d$ can be identified by Algorithm~\ref{alg:order-d} up to gauge invariance in \eqref{eq:gauge-invariance}.
\end{theorem}

\begin{algorithm}[ht]
    \caption{Blockwise Simultaneous Diagonalization Method for Tensor Ring Decomposition (\blostr)}
    \label{alg:order-d}
    \KwIn{Tensor $\mbcT$ observable at entries $\Delta$ in \eqref{eq:Delta}, TR-rank $r$, index $\bgamma\in[n_1]\times\cdots\times[n_d]$}
    \KwOut{TR-cores $\hat{\mbcQ}_k\in\bbC^{n_k\times r\times r}$, $k\in[d]$}
    Choose two index pairs $(\balpha,\bbeta)$ and $(\balpha',\bbeta')$ in $([n_2]\times\cdots\times[n_{d-1}])^2$ satisfying \eqref{eq:alpha-beta}
    \textit{(Practical: try multiple index pairs and keep one with well-separated spectra in Step 2.)}\;
    
    Calculate the eigendecomposition of $\bfT(:,\balpha,\Gamma_{\balpha})\bfT(:,\bbeta,\Gamma_{\bbeta})^\dagger$ and $\bfT(:,\balpha',\Gamma_{\balpha'})\bfT(:,\bbeta',\Gamma_{\bbeta'})^\dagger$. Compute eigenvectors corresponding to the nonzero spectrum and stack them as $\bfE,\bfE'\in\bbC^{n_1\times r^2}$\;
    
    Partition columns of $\bfE$ (resp.\ $\bfE'$) into $r$ blocks of size $r$ so that each block contains eigenvectors corresponding to the $r$ distinct nonzero eigenvalues (up to numerical error) in a fixed order\;
    
    Set $\bfF=\bfE^\dagger\bfE'$ and define the $(j,k)$-th \(r\times r\) block by $\bfF^{(j,k)}=\bfF_{[0:(r{-}1)]\cdot r+j,\,[0:(r{-}1)]\cdot r+k}$ for $j,k\in[r]$\;
    
    Set $\hat{\bfK}_1=\bfI_r$. For $\ell=2,\dots,r$, set $\hat{\bfK}_\ell=\bfF^{(1,1)}\big(\bfF^{(\ell,1)}\big)^{-1}$ and let $\hat{\bfK}=\mathrm{blockdiagonal}(\hat{\bfK}_1,\dots,\hat{\bfK}_r)$\;
    
    Recover the first core unfolding:
    $\hat{\bfQ}_{1\langle1\rangle}=\bfE\big(\Pi(\hat{\bfK})\big)^{-1}$,
    then reshape to obtain $\hat{\mbcQ}_1$\;
    
    \For{$k=2$ \KwTo $d$}{
      Compute
      \[
      \hat{\bfQ}_{k\left<1\right>}=\overleftarrow{\bfT}^{k-1}(:,\overleftarrow{\bgamma}_{\rm mid}^{k-1},\Gamma_{k-1})\left(\left(\bfI_r\otimes(\hat{\bfQ}_1^{(\gamma_1)}\cdots\hat{\bfQ}_{k-2}^{(\gamma_{k-2})})\right)\hat{\bfQ}_{(k-1)[1]}^*(:,\Gamma_{k-1})\right)^\dagger,
      \]
      with the convention that the product in parentheses is $\bfI_r$ when $k=2$\;
      
      Reshape to obtain $\hat{\mbcQ}_k$
    }
\end{algorithm}

\begin{figure}[htbp]
    \centering
    \includegraphics[width=\linewidth]{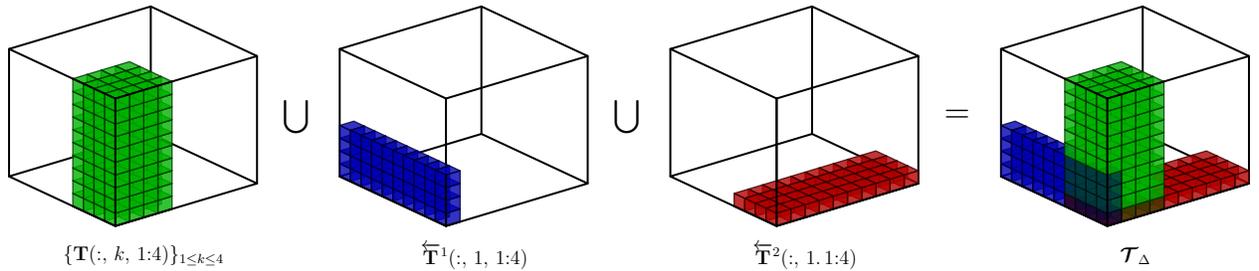}
    \caption{Illustrative example of a possible choice of $\Delta$ in a three-way tensor. Each small cell represents an entry. We set $n_1=n_2=n_3=10$, $r=2$, and all $\Gamma_{\cdot}$'s are $\{1,2,3,4\}$.}
    \label{fig:three-tensor-delta}
\end{figure}

Theorem~\ref{thm:order-d} establishes two key advantages of Algorithm~\ref{alg:order-d}. 
First, it recovers the exact TR decomposition in a finite number of steps. 
Second, it requires only $O((n_1+\cdots + n_d)r^2)$ tensor entries to recover the TR cores, which is optimal in terms of sampling size. In particular, for any $\balpha,\bbeta,\balpha',\bbeta'$ satisfying \eqref{eq:alpha-beta} and $\bgamma\in[n_1]\times\cdots\times[n_d]$, one only needs to observe the following $O((n_1+\cdots + n_d)r^2 -dr^2)$ entries 
\begin{align}
    \label{eq:Delta}
    \mbcT_{\Delta}\;\equiv\; \bigl\{\bfT(:,\balpha,\Gamma_{\balpha}), \bfT(:,\bbeta,\Gamma_{\bbeta}) \bigr\}\,\bigcup\,\bigl\{ \bfT(:,\balpha',\Gamma_{\balpha'}),\,\bfT(:,\bbeta',\Gamma_{\bbeta'})\bigr\}\, \bigcup\, \bigl\{\overleftarrow{\bfT}^k(:, \overleftarrow{\bgamma}_{\rm mid}^k, \Gamma_k)\bigr\}_{k\in[d]},
\end{align}
where $\Delta\subset[n_1]\times\cdots\times[n_d]$ denotes the set of corresponding indices. Figure~\ref{fig:three-tensor-delta} demonstrates a possible choice of $\mbcT_\Delta$ for a three-way tensor. This sampling efficiency is particularly valuable for high-order or high-dimensional tensors. Note that the TR decomposition has $O((n_1+\cdots +n_d)r^2)$ free parameters up to gauge invariance. This further shows {\bf \blostr achieves the optimal sample complexity.}

\begin{remark}[Computational complexity]
    The dominant cost of Algorithm~\ref{alg:order-d} is the eigendecomposition step in Line~2, which operates on an $n_1 \times n_1$ matrix of rank $r^2$, costing $O(n_1 r^4)$. The sequential core recovery in Lines~7--9 costs $O(\sum_k n_k r^4)$ in total, giving an overall complexity of $O\left((\sum_k n_k) r^4\right)$, which is polynomial in all problem parameters.
\end{remark}

\begin{remark}
    The condition $n_k \geq r^2$ for all $k\in[d]$ in Theorem~\ref{thm:order-d} can be relaxed to $n_k\geq r^2$ for some $k\in[d]$, at the expense of requiring more observed entries of the tensor. The refined algorithm is presented in Section~\ref{sec:algorithm}.
\end{remark}

\begin{remark}
    For simplicity, Theorem~\ref{thm:order-d} assumes all TR ranks are identical. A slight modification of Algorithm~\ref{alg:order-d} allows one TR rank to differ from the others. The more general case in which all ranks are different is a more challenging problem and remains open; we leave it as future work.
\end{remark}

\section{BLOSTR for Symmetric Tensor-Ring Decomposition}\label{sec:symmetric}

We next consider the \emph{symmetric} TR setting, where $\mbcT\in\bbC^{n^{\times d}}$, and for all $\balpha=(\alpha_1,\alpha_2,\ldots,\alpha_d)\in[n]^d$,
\begin{align}\label{eq:symmetric-n}
  T(\alpha_1,\alpha_2,\ldots,\alpha_d)
  \;=\;
  \mathrm{tr}\!\left\{\bfQ^{(\alpha_1)}\bfQ^{(\alpha_2)}\cdots\bfQ^{(\alpha_d)}\right\}.
\end{align}

Symmetric TR decomposition in \eqref{eq:symmetric-n} is motivated by applications where invariance or homogeneity plays a central role. In statistics and machine learning, it is well suited for exchangeable or permutation-invariant settings (e.g., graphs or i.i.d. random variables), where tying all cores yields a more structured and interpretable representation \cite{chen2023learning,comon1996decomposition,zhou2006learning}. In quantum many-body physics, it naturally models translationally invariant systems such as bosonic wavefunctions, which are symmetric under the exchange of identical particles \cite{marconi2025symmetric}. From a computational perspective, symmetry also leads to substantial efficiency gains: while a general TR with $d$ modes requires $O(dnr^2)$ parameters, the symmetric TR reduces this to only $O(nr^2)$, making it particularly attractive for high-order problems.

Suppose $\{\hat{\mbcQ}_k\}_{k=1}^d$ is any (not necessarily symmetric) TR decomposition of $\mbcT$ obtained by Algorithm~\ref{alg:order-d}.
We seek a single core $\tilde{\mbcQ}$ that reproduces all entries via \eqref{eq:symmetric-n}.
The following consistency statement is immediate.

\begin{lemma}\label{lemma:hat-q-tilde}
Let $\{\hat{\mbcQ}_k\}_{k=1}^d$ satisfy the TR decomposition \eqref{eq:order4TR}. If there exists $\tilde{\mbcQ}\in\bbC^{n\times r\times r}$ such that
\begin{align}\label{eq:tilde.q}
  \hat{\bfQ}_1^{(\alpha_1)}\hat{\bfQ}_2^{(\alpha_2)}\cdots\hat{\bfQ}_d^{(\alpha_d)}
  \;=\;
  \tilde{\bfQ}^{(\alpha_1)}\tilde{\bfQ}^{(\alpha_2)}\cdots\tilde{\bfQ}^{(\alpha_d)},
  \quad \forall\,(\alpha_1,\ldots,\alpha_d)\in[n]^d,
\end{align}
then $\tilde{\mbcQ}$ is a symmetric TR core for $\mbcT$ since \eqref{eq:symmetric-n} holds.
\end{lemma}

Note that if we know $\tilde{\bfQ}^{(1)}$, then the rest of the slices can be obtained by letting $\balpha=(1,\dots,1,\alpha_d)$ in \eqref{eq:tilde.q}. Therefore, \eqref{eq:tilde.q} defines a system of polynomial equations with $n^d$ equations in $nr^2$ variables. A natural first approach is to consider Gr\"obner bases, which provide a canonical generating set for polynomial systems. However, computing a Gr\"obner basis can require doubly exponential time in the number of variables in the worst case \cite{mayr1982complexity}. Although modern algorithms such as Faug\`ere’s F4 and F5 have demonstrated practical efficiency on structured instances \cite{cox2015ideals}, the systems arising from \eqref{eq:tilde.q} typically exhibit little exploitable algebraic structure, since both the tensor $\mbcT$ and the estimator $\hat{\mbcQ}$ are random. As a consequence, the construction of the system is already computationally burdensome, and the absence of structure further exacerbates the difficulty of computing the corresponding Gr\"obner basis.

Alternatively, we consider searching over the possible space of solutions based on several observations. For the all-equal tuples $(\alpha,\ldots,\alpha)$, \eqref{eq:tilde.q} reduces to $\hat{\bfQ}_1^{(\alpha)}\cdots\hat{\bfQ}_d^{(\alpha)}\;=\; \big(\tilde{\bfQ}^{(\alpha)}\big)^{d}$. This yields a constructive route via diagonalization as described by the following lemma.

\begin{lemma}
    \label{lemma:q}
    For any $\alpha\in[n]$, assume that $\hat{\bfQ}_1^{(\alpha)}\hat{\bfQ}_2^{(\alpha)}\cdots\hat{\bfQ}_d^{(\alpha)}$ has an eigendecomposition as $\bfY_\alpha\boldsymbol{\Lambda}_\alpha\bfY_\alpha^{-1}$, where $\bfY_\alpha\in\text{GL}(r,\bbC)$ and $\boldsymbol{\Lambda}_\alpha=\text{diag}(\lambda_{\alpha,1},\dots,\lambda_{\alpha,r})$. Then $\tilde{\bfQ}^{(\alpha)}$ is in the form of $\bfY_\alpha\boldsymbol{\Omega}_\alpha\bfY_\alpha^{-1}$, where $\boldsymbol{\Omega}_\alpha=\text{diag}(\omega_{\alpha,1},\dots,\omega_{\alpha,r})$ and $\omega_{\alpha,\ell}$ is a $d$-th root of $\lambda_{\alpha,\ell}$.
\end{lemma}

The existence of $\tilde{\mbcQ}$ is guaranteed by the existence of $\mbcQ$ in \eqref{eq:symmetric-n}. The following lemma indicates that the solution is not unique.

\begin{lemma}\label{lemma:roots}
    Assume that $\boldsymbol{\Lambda}_1=\text{diag}(\lambda_{1,1},\dots,\lambda_{1,r})$ and $\lambda_{1,t}=\ell_te^{i\theta_t}$, where $\ell_t\geq0$ and $\theta_t\in[0,2\pi)$. If $\bar{\bfQ}^{(1)}=\bfY_1\boldsymbol{\Omega}_1\bfY_1^{-1}$ is a solution for \eqref{eq:tilde.q}, where $\boldsymbol{\Omega}_1=\text{diag}(\omega_{1,1},\dots,\omega_{1,r})$ and $\omega_{1,t}=\ell_t^{1/d}e^{i(\theta_t+2k_t\pi)/d}$ for some $k_t\in\{0,1,\dots,d-1\}$, then the matrices $\boldsymbol{\Omega}_1'=\text{diag}(\omega_{1,1}',\dots,\omega_{1,r}')$ with $\omega_{1,t}'=\ell_t^{1/d}e^{i(\theta_t+2(k_t+1)\pi)/d}$ is also a solution to \eqref{eq:tilde.q}.
\end{lemma}

Combining Lemma~\ref{lemma:q} and Lemma~\ref{lemma:roots} allows us to restrict the (equivalent) solution of \eqref{eq:tilde.q} to the set
\begin{align}
    \label{eq:solution-set}
    \left\{\boldsymbol{\Omega}=\text{diag}(\omega_{1,1},\dots,\omega_{1,r})\ |\ \omega_{1,1}=\ell_1^{1/d} e^{i\theta_1/d},\ \omega_{1,t}^d=\lambda_{1,t}\ \text{for}\ 2\leq t\leq r\right\}.
\end{align}
We formalize our discussion in the following theorem.

\begin{theorem}[Blockwise Simultaneous Diagonalization Method for Symmetric Tensor Ring Decomposition]\label{thm:symmetric}
    Suppose $\mbcT\in\bbC^{n^{\times d}}$ satisfies \eqref{eq:symmetric-n} with $n\geq r^2$. Assume that the elements of $\mbcQ\in\mathbb{C}^{n\times r\times r}$ are randomly drawn from a measure $\mu$ that is absolutely continuous with respect to the standard Euclidean measure. Then with probability one, $\mbcQ$ can be identified by Algorithm~\ref{alg:symmetric} up to gauge invariance in \eqref{eq:gauge-invariance}.
\end{theorem}

\begin{algorithm}[ht]
    \caption{\texttt{BLOSTR-S} for Symmetric Tensor Ring Decomposition}
    \label{alg:symmetric}
    \KwIn{Symmetric tensor $\mbcT$ observable at entries $\Delta_s$ in \eqref{eq:Delta-s}, TR-rank $r$}
    \KwOut{TR-core $\tilde{\mbcQ}\in\bbC^{n\times r\times r}$}
    Implement Algorithm~\ref{alg:order-d} to obtain $\{\hat{\mbcQ}_k\}_{k\in[d]}$\;
    
    Obtain the eigendecomposition $\hat{\bfQ}_1^{(1)}\hat{\bfQ}_d^{(1)}=\bfY_1\boldsymbol{\Lambda}_1\bfY_1^{-1}$, with $\boldsymbol{\Lambda}_1=\mathrm{diag}(\lambda_{1,1},\ldots,\lambda_{1,r})$\;
    
    Obtain $(\ell_1,\dots,\ell_r)$ and $(\theta_1,\dots,\theta_r)\in[0,2\pi)^r$ such that $\lambda_{1,k}=\ell_ke^{i\theta_k}$\;
    
    \Repeat{
        $\bar{\mbcQ}$ satisfies \eqref{eq:tilde.q}
    }
    {
        Choose $(k_2,\dots,k_r)\in\{0,1,\dots,d-1\}^{r-1}$\;
        
        let $\boldsymbol{\Omega}=\mathrm{diag}\left(\ell_1^{1/d}e^{i\theta_1/d},\ 
                \ell_2^{1/d}e^{i(\theta_2+2k_2\pi)/d},\ \dots,\ 
                \ell_r^{1/d}e^{i(\theta_r+2k_r\pi)/d}\right)$\;
        
        $\bar{\bfQ}^{(1)}\gets \bfY_1\boldsymbol{\Omega}\bfY_1^{-1}$\;
        
        $\bar{\bfQ}^{(j)}\gets (\bar{\bfQ}^{(1)})^{-(d-1)}\hat{\bfQ}_1^{(1)}\hat{\bfQ}_d^{(j)}$ for $j\in\{2,\ldots,n\}$\;
    }
    \Return{$\tilde{\mbcQ}\gets \bar{\mbcQ}$}
\end{algorithm}

Algorithm~\ref{alg:symmetric} does not require any additional observations compared to Algorithm~\ref{alg:order-d}. In fact, compared with \eqref{eq:Delta}, the set of entries of $\mbcT$ to be observed reduces to $O(nr^2)$:
\begin{align}
    \label{eq:Delta-s}
    \mbcT_{\Delta_s}\;\equiv\; \bigl\{\bfT(:,\balpha,\Gamma_{\balpha}),\,\bfT(:,\bbeta,\Gamma_{\bbeta}) \bigr\}\,\bigcup\,\bigl\{\bfT(:,\balpha',\Gamma_{\balpha'}),\,\bfT(:,\bbeta',\Gamma_{\bbeta'})\bigr\}\,\bigcup\, \bigl\{\bfT(:,\bgamma,\Gamma_{\bgamma})\bigr\}
\end{align}
for some $\balpha,\balpha',\bbeta,\bbeta'$ satisfying \eqref{eq:alpha-beta}, and $\bgamma=m\cdot\boldsymbol{1}_{d-2}$ for $m\leq n$. This aligns with the degrees of freedom of a symmetric TR decomposition up to gauge invariance, namely $O(nr^2)$. Hence, \texttt{BLOSTR-S}  attains the optimal sample complexity.

\section{TR Decomposition in the Presence of Perturbations}\label{sec:robust}

In this section, we study the TR decomposition in a practical setting where the tensor is observed with perturbations. 
Since our analysis does not rely on a specific assumption about the perturbations, we employ a pragmatic procedure to handle them, and we demonstrate in Section~\ref{sec:simulation} that our proposed method performs well empirically. We leave a rigorous quantitative analysis to future work.
In particular, we aim to identify TR components $\mbcQ_k\in\bbC^{n_k\times r\times r}$, $k\in[d]$ that approximately decompose $\mbcT$:
\begin{align}\label{eq:tr-noise}
    T(\alpha_1,\alpha_2,\dots,\alpha_d) = \text{tr}\left\{\bfQ_1^{(\alpha_1)}\bfQ_2^{(\alpha_2)}\cdots\bfQ_d^{(\alpha_d)}\right\} + W(\alpha_1, \alpha_2, \ldots, \alpha_d), \quad \alpha_k \in [n_k].
\end{align}
Here, $W(\alpha_1, \alpha_2, \ldots, \alpha_d)$ represents the perturbation. We estimate a rank-$r$ TR model by approximately solving
\begin{equation}\label{eq:opt}
\begin{split}
    (\hat{\mbcQ}_1,\hat{\mbcQ}_2,\dots,\hat{\mbcQ}_d) = & \argmin_{\mbcQ_1,\mbcQ_2,\dots,\mbcQ_d}\left\|\left\{\mbcT-\mfR(\mbcQ_1,\mbcQ_2,\dots,\mbcQ_d)\right\}_{\Delta}\right\|_F^2\\
    = & \argmin_{\mbcQ_1,\mbcQ_2,\dots,\mbcQ_d} \sum_{(j_1,\ldots, j_d)\in \Delta}|\mbcT_{j_1,\ldots, j_d} - \mfR(\mbcQ_1,\mbcQ_2,\dots,\mbcQ_d)_{j_1,\ldots, j_d}|^2,
\end{split}
\end{equation}
where $\Delta$ denotes the indices of the observed entries of $\mbcT$ defined in \eqref{eq:Delta}. 

We describe the step-wise procedure below. First, we adapt the procedure of \blostr to this setting. Recall that for any tensor $\mbcT\in\bbC^{n_1\times \cdots \times n_d}$ and $\balpha\in[n_2]\times\cdots\times [n_{d-1}]$, $\bfT( :, \balpha, \Gamma_{\balpha})$ is the $n_1 \times r^2$ matrix obtained by fixing the mode-$2$ through mode-$(d-1)$ indices to $\balpha$ and restricting the mode-$d$ index to $\Gamma_{\balpha}\subseteq [n_d]$. As in Section~\ref{sec:procedure}, choose indices
$\balpha,\bbeta\in[n_2]\times\cdots\times[n_{d-1}]$, $\balpha\neq\bbeta$, and form
\[
  \bfM \;=\; \bfT(:,\balpha,\Gamma_{\balpha})\big(\bfT(:,\bbeta,\Gamma_{\bbeta})\big)^\dagger.
\]
Under perturbations, the $r$ eigenvalues with multiplicity $r$ of $\bfM$ in the noiseless case generally separate into $r^2$ eigenvalues that form $r$ clusters. To recover these groups in practice, we cluster the $r^2$ eigenvalues of \(\bfM\) with the largest modulus into $r$ equal-size clusters in the complex plane (embedding $x_j=(\Re(\lambda_j),\Im(\lambda_j))\in\mathbb{R}^2$). We perform this grouping using constrained $k$-means~\cite{bradley2000constrained}:
\begin{equation}\label{eq:kmeans}
  \min_{\{z_{j,k}\},\,\{c_k\}}
  \sum_{j=1}^{r^2}\sum_{k=1}^{r} z_{j,k}\,\|x_j-c_k\|_2^2
  \ \ \text{s.t.}\ \
  \sum_{k=1}^{r} z_{j,k}=1,\;
  \sum_{j=1}^{r^2} z_{j,k}=r,\;
  z_{j,k}\in\{0,1\}.
\end{equation}
Since constrained clustering is a combinatorial problem, we use a standard heuristic solver\footnote{We use SVD-based pseudoinverses and an eigenvalue clustering tolerance; an implementation option is \texttt{k-means-constrained}~\cite{Levy-Kramer_k-means-constrained_2018}.}; in our experiments it is stable when the spectral separation between groups exceeds the noise level. To mitigate sensitivity to initialization, we run the clustering with multiple random restarts and retain the solution with the smallest objective value in \eqref{eq:kmeans}. The resulting block structure yields an estimate of the eigenvector matrix $\bfE$ (and, analogously, $\bfE'$ from a second pair $(\balpha',\bbeta')$), from which we construct initial cores by the same block-recovery steps as in Section~\ref{sec:procedure}.

Next, given an initial set $\{\mbcQ_k\}_{k=1}^d$, we refine the estimate by alternating least squares (ALS) following~\cite{zhao2016tensor}.
For $k\in[d]$, define the contracted environment tensor $\mbcQ^{\neq k}\in\bbC^{(\prod_{j\ne k} n_j)\times r\times r}$ slice-wise by
\[
  \bfQ^{\neq k(\overline{\alpha_{k+1}\cdots \alpha_d\,\alpha_1\cdots \alpha_{k-1}})}
  \;=\;
  \Bigl(\prod_{\ell=k+1}^{d}\bfQ_\ell^{(\alpha_\ell)}\Bigr)
  \Bigl(\prod_{m=1}^{k-1}\bfQ_m^{(\alpha_m)}\Bigr),
\]
and update $\mbcQ_k$ by solving the linear least-squares problem
\begin{align}
    \label{eq:als-q}
    \min_{\bfQ_{k[1]}} \ \bigl\| \bigl(\bfT_{[k]} - \bfQ_{k[1]}\,(\bfQ^{\neq k}_{\langle 1\rangle})^* \bigr)_{\Delta}\bigr\|_F^2,
\end{align}
where $\bfT_{[k]}$, $\bfQ_{k[1]}$ and $\bfQ_{\langle 1\rangle}^{\neq k}$ are the matricization of $\mbcT$, $\mbcQ_k$ and $\mbcQ^{\neq k}$ defined in \eqref{eq:matricization}, respectively; and where (with a slight abuse of notation) $\Delta$ denotes the corresponding index set after matricization. For $k\in[d]$ and $\ell\in[n_k]$, denote by $\bfq_{k,\ell}$ the $\ell$th row of $\bfQ_{k[1]}$. Let 
\begin{align*}
    \mcJ_{k,\ell}\equiv\left\{\overline{j_{k+1} \cdots j_d j_1 \cdots j_{k-1}} \mid (j_1,\dots,j_{k-1},\ell,j_{k+1},\dots,j_d)\in\Delta \right\}.
\end{align*} 
Equation \eqref{eq:als-q} is equivalent to the following least squares problem:
\begin{align}
    \label{eq:masked-ls}
    \min_{\bfq_{k,\ell} \in\bbC^{r^2}} \|\bfT_{[k]}(\ell,\mcJ_{k,\ell}) - \bfq_{k,\ell}^\top \bigl(\bfQ_{\langle1\rangle}^{\neq k}(\mcJ_{k,\ell},:)\bigr)^*\|_2^2,\quad \ell\in[n_k],\quad \mcJ_{k,\ell}\neq\emptyset.
\end{align} 
Closed-form updates use normal equations or QR/SVD solvers. 
We iterate over $k=1,\ldots,d$ until convergence (relative decrease below a tolerance or a maximum number of sweeps). The overall procedure is summarized in Algorithm~\ref{alg:noise}. 

\begin{algorithm}[ht]
    \caption{\blostr for Robust TR Decomposition}
    \label{alg:noise}
    \KwIn{Tensor $\mbcT\in\bbC^{n_1\times n_2\times\cdots\times n_d}$ observed at entries $\Delta$ in \eqref{eq:Delta}, TR-rank $r$}
    \KwOut{$\hat{\mbcQ}_k\in\bbC^{n_k\times r\times r}$, $k\in[d]$, solution to \eqref{eq:opt}}
    
    Choose two index pairs $(\balpha,\bbeta)$ and $(\balpha',\bbeta')$ satisfying \eqref{eq:alpha-beta}\;
    
    Obtain SVD $\bfT(:,\balpha,\Gamma_{\balpha})=\bfU_{\balpha}\boldsymbol{\Sigma}_{\balpha}\bfV_{\balpha}^*$ and $\bfT(:,\bbeta,\Gamma_{\bbeta})=\bfU_{\bbeta}\boldsymbol{\Sigma}_{\bbeta}\bfV_{\bbeta}^*$, where $\boldsymbol{\Sigma}_{\balpha},\boldsymbol{\Sigma}_{\bbeta}\in\bbC^{r^2\times r^2}$ are the top $r^2$ singular values of $\bfT(:,\balpha,\Gamma_{\balpha})$ and $\bfT(:,\bbeta,\Gamma_{\bbeta})$, respectively, and $\bfU_{\balpha},\bfU_{\bbeta}\in\bbC^{n_1\times r^2}$, $\bfV_{\balpha},\bfV_{\bbeta}\in\bbC^{r^2\times r^2}$ are the corresponding column and row vectors\;
    
    Obtain the eigendecomposition of $\bfU_{\balpha}\boldsymbol{\Sigma}_{\balpha}\bfV_{\balpha}^*\bfV_{\bbeta}\boldsymbol{\Sigma}_{\bbeta}^{-1}\bfU_{\bbeta}^*$\;
    
    Apply constrained $k$-means clustering \eqref{eq:kmeans} on the $r^2$ eigenvalues with the largest modulus to get $r$ clusters of size $r$. Rearrange the corresponding eigenvectors to obtain the new matrices $\bfE\in\bbC^{n_1\times r^2}$, such that the eigenvalue $\lambda_{(j-1)r+k}$ is in the $k$-th cluster, $j,k\in[r]$\;
    
    Repeat the previous three steps above for $\bfT(:,\balpha',\Gamma_{\balpha'})$ and $\bfT(:,\bbeta',\Gamma_{\bbeta'})$ to obtain $\bfE'\in\bbC^{n_1\times r^2}$\;
    
    Follow lines 3--9 in Algorithm~\ref{alg:order-d} with $\bfE$ and $\bfE'$ to obtain $\hat{\mbcQ}_{k,0}\in\bbC^{n_k\times r\times r}$, $k\in[d]$\; 
    
    \tcp{The second subscript in $\hat{\mbcQ}_{k,0}$ indicates the number of iteration, $t=0$}
    
    Set $t \gets 0$.

\While{ $\|\{\mbcT-\mfR(\hat{\mbcQ}_{1,t},\dots,\hat{\mbcQ}_{d,t})\}_\Delta\|_F \ge \epsilon$ and $t < t_{\max}$ }{
  Sequentially update $\hat{\mbcQ}_{k,t}$ by \eqref{eq:masked-ls} for $k\in[d]$ to obtain $\hat{\mbcQ}_{k,t+1}$\;
  Set $t \gets t+1$\;
}
    \Return{$\hat{\mbcQ}_k \gets \hat{\mbcQ}_{k,t}$, $k\in[d]$}
\end{algorithm}

\begin{remark}[TR Decomposition in the Presence of Perturbations]
Furthermore, for a perturbed observation of a tensor that admits a symmetric decomposition in the sense of \eqref{eq:symmetric-n}, one may first compute $\{\hat{\mbcQ}_k\}_{k\in[d]}$ using Algorithm~\ref{alg:noise} with the input of entries $\mbcT_{\Delta_s}$. Subsequently, $\hat{\mbcQ}$ can be obtained by selecting $\boldsymbol{\Omega}$ from the solution set \eqref{eq:solution-set} to minimize the Frobenius norm discrepancy:
\[
    \sum_{\balpha\in[n]^d}\left\| \hat{\bfQ}_1^{(\alpha_1)}\hat{\bfQ}_2^{(\alpha_2)}\cdots\hat{\bfQ}_d^{(\alpha_d)} - \hat{\bfQ}^{(\alpha_1)}\hat{\bfQ}^{(\alpha_2)}\cdots\hat{\bfQ}^{(\alpha_d)}\right\|_F^2.
\]
\end{remark}

\section{Numerical Studies}\label{sec:simulation}

In this section, we report numerical experiments evaluating the performance of the proposed \blostr algorithms. Due to space constraints, additional results are provided in Appendix~\ref{sec:addsimulation}. 

We first assess Algorithm~\ref{alg:order-d} across tensor orders $d\in\{5,6,7\}$ and a range of TR rank settings. The ground-truth TR cores $\{\mbcQ_k\}_{k\in[d]}$ have entries drawn i.i.d.\ from $\mathcal{N}(0,10^2)$. As summarized in Table~\ref{tab:exact}, the reconstruction error
$\|\mbcT-\hat{\mbcT}\|_F$, with $\hat{\mbcT}$ produced by Algorithm~\ref{alg:order-d} (or, alternatively, Algorithm~\ref{alg:refined}), is negligible up to machine precision, indicating exact recovery of the TR decomposition.

\begin{table}[htbp]
    \centering
    \caption{Reconstruction errors of exact TR decomposition under various dimensions and TR rank settings. The entries of the TR cores, $\{\mbcQ_k\}_{k\in[d]}$, are drawn independently from $\mcN(0, 10^2)$.}
    \begin{tabular}{llcc}
        \toprule
        $\bfn=(n_1,n_2,\dots,n_d)$ & $r$ & $\|\mbcT-\hat{\mbcT}\|_F$ & $\|\mbcT-\hat{\mbcT}\|_F/\|\mbcT\|_F$\\
        \midrule
        (12, 5, 6, 7, 10) & 3 & 1.09e-03 & 5.73e-12 \\
        (10, 10, 10, 10, 10) & 2 & 8.89e-05 & 5.70e-13 \\
        (20, 20, 20, 20, 20) & 2 & 4.49e-05 & 6.12e-14 \\
        (20, 20, 20, 20, 20) & 4 & 2.93e-02 & 5.72e-12 \\
        (10, 10, 10, 10, 10, 10) & 2 & 2.96e-03 & 3.83e-13 \\
        (10, 10, 10, 10, 10, 10, 10) & 2 & 1.54e+00 & 4.81e-12 \\
        \bottomrule
    \end{tabular}
    \label{tab:exact}
\end{table}

\begin{figure}[ht]
    \centering
    \includegraphics[width=\linewidth]{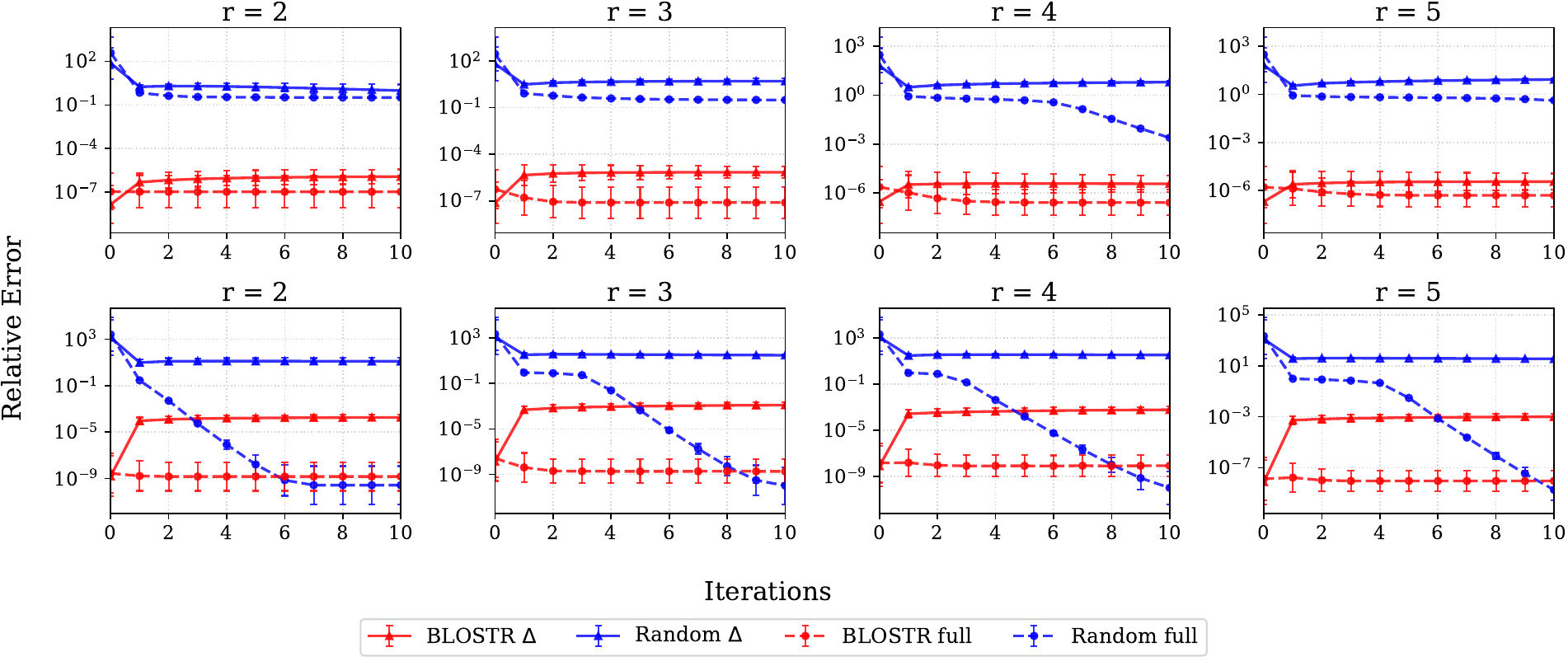}
    \caption{Average relative error (with standard deviation in log scale) of \blostr (Algorithm~\ref{alg:noise}) and randomly initialized ALS over varying numbers of iterations. Full: access to full entries of tensor $\mbcT$ at the stage of ALS; $\Delta$: ALS with entries limited to $\Delta$ in \eqref{eq:Delta}. The entries of $\{\mbcQ_k\}_{k\in[d]}$ are independently drawn from $\mcN(0, 10^2)$, and the variance of the added Gaussian noise is set to 1. The tensor $\mbcT$ has dimension $30^{\times 3}$ in the top row and $30^{\times 4}$ in the bottom row.}
    \label{fig:err-vs-iter}
\end{figure}

We benchmark Algorithm \ref{alg:noise} against a standard TR-ALS baseline with random initialization over 100 independent trials. For each trial, we record the relative reconstruction error
$\|\mbcT - \hat{\mbcT}\|_F / \|\mathcal{T}\|_F$
over 10 ALS iterations on tensors of size $30^{\times 3}$ and $30^{\times 4}$ with TR rank $r$, corrupted by additive Gaussian noise of variance~$1$.
As shown in Figure~\ref{fig:err-vs-iter}, \blostr attains substantially lower reconstruction error in the early iterations and converges in fewer iterations, which can be attributed to its structured initialization.

\begin{figure}[ht]
    \centering
    \includegraphics[width=\linewidth]{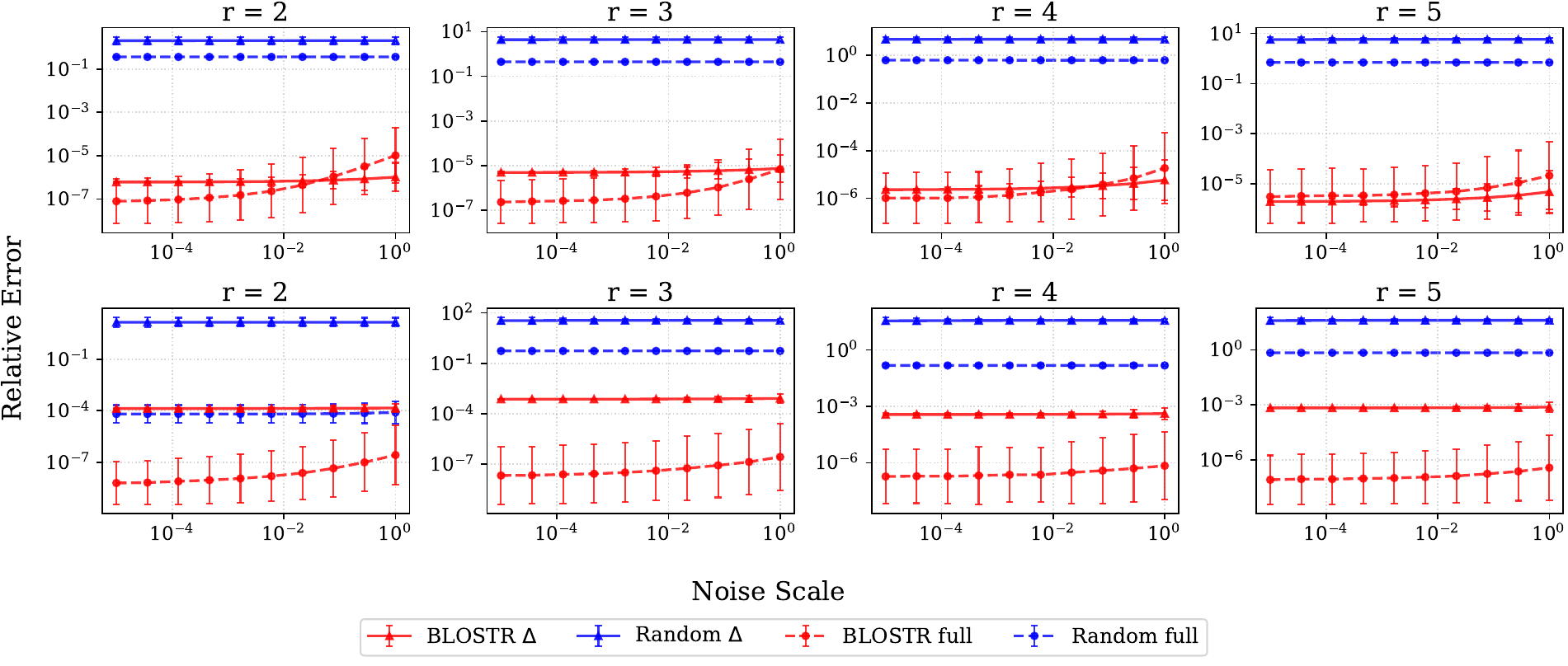}
    \caption{Average relative error (with standard deviation in log scale) of \blostr (Algorithm~\ref{alg:noise}) and randomly initialized ALS over increasing noise scales with 3 iterations. Full: access to full entries of tensor $\mbcT$ at the stage of ALS; $\Delta$: ALS with entries limited to $\Delta$ in \eqref{eq:Delta}. The entries of $\{\mbcQ_k\}_{k\in[d]}$ are independently drawn from $\mcN(0, 10^2)$. The tensor $\mbcT$ has dimension $30^{\times 3}$ in the top row and $30^{\times 4}$ in the bottom row.}
    \label{fig:err-vs-ns-t3}
\end{figure}

We further examine the robustness of the proposed algorithm under increasing noise scales. Figure~\ref{fig:err-vs-ns-t3} shows that \blostr consistently achieves low reconstruction error within only three iterations across all noise levels, with mean relative errors notably smaller than those obtained by random initialization.

\begin{figure}[ht]
    \centering
    \includegraphics[width=\linewidth]{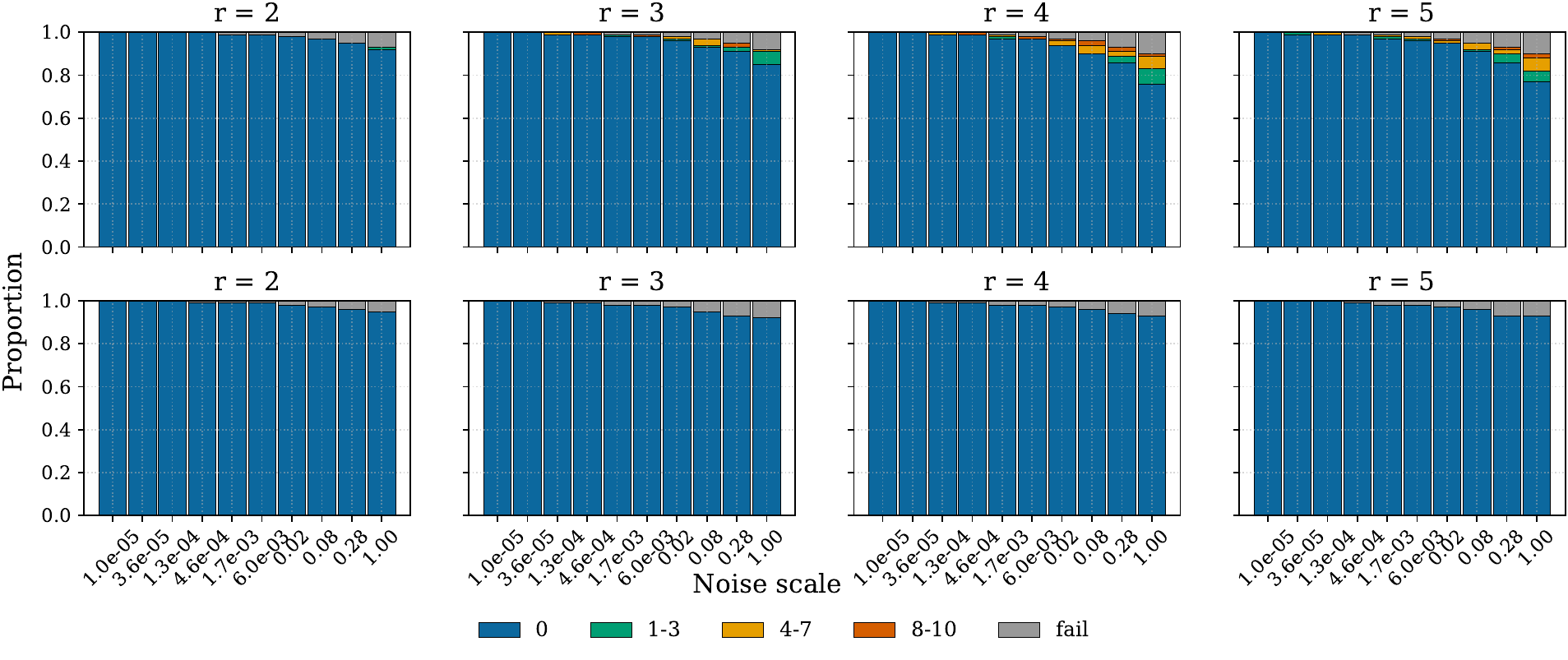}
    \caption{Proportion of successful recoveries (relative error below $10^{-5}$) within at most 10 iterations, averaged over 100 independent trials, using \blostr (Algorithm \ref{alg:noise}) under varying noise scales $\sigma_n$. Each tensor core $\{\mbcQ_k\}_{k\in[d]}$ has entries independently drawn from $\mcN(0, 10^2)$. The underlying tensor $\mbcT$ has dimension $30^{\times 3}$ in the top row and $30^{\times 4}$ in the bottom row.}
    \label{fig:t2tol-1e-5}
\end{figure}

Finally, Figure~\ref{fig:t2tol-1e-5} summarizes the proportion of successful recoveries (relative error below $10^{-5}$) within at most 10 iterations across different noise scales. The results indicate that \blostr reliably recovers the underlying tensor when the noise is small and maintains a high probability of success even at moderate noise levels.

Overall, these results demonstrate that \blostr achieves reliable and efficient tensor ring decomposition.

\section{Connections to Matrix Product States and Moment Tensors}\label{sec:mps}

In this section, we highlight connections between TR decomposition and two other well-studied problems, one from quantum information and the other from high-dimensional statistics. The first connection is to \emph{matrix product state (MPS) tomography}, which we elaborate upon in Section~\ref{sec:main_mps}. The second connection is to \emph{learning polynomial transformations}, which we elaborate upon in Section~\ref{sec:main_moment}. These two perspectives illustrate how TR serves both as a tensor network structure and as a statistical tool, while leaving open further connections to be explored in other domains.

\subsection{Matrix Product States}
\label{sec:main_mps}

In quantum many-body physics, the state of a physical system given by $d$ particles is described by a unit vector in a Hilbert space of dimension exponential in $d$, e.g., $\mathbb{C}^{n_1}\otimes \cdots \otimes \mathbb{C}^{n_d}$. Whereas a generic state thus requires exponentially many parameters to describe, there are a variety of ways to parameterize the set of ``physically relevant'' states using only a polynomial number of parameters, one prominent example being \emph{matrix product states}. Concretely, a matrix product state, denoted in bra-ket notation by $\ket{\psi}$, is given by a complex unit vector in the aforementioned space, whose entries are given by the entries of a tensor $\mbcT$ with bounded TR rank, which in this literature is called the \emph{bond dimension}.

\begin{definition}[Matrix product states]
    A unit vector $\ket{\psi}$ is an MPS with \emph{bond dimension $r$} if its entries are given by
    \begin{equation}
        \label{eq:MPS}
        \ket{\psi} = \sum_{\alpha_1,\ldots,\alpha_d} \mathrm{tr}\!\left\{\bfQ^{(\alpha_1)}_1\bfQ^{(\alpha_2)}_2\cdots \bfQ^{(\alpha_d)}_d\right\}\ket{\alpha_1,\ldots,\alpha_d},
    \end{equation}
    where each $\bfQ^{(\alpha_k)}_k$ is an $r\times r$ matrix, and $\ket{\alpha_1,\ldots,\alpha_d}$ denotes the standard basis vector of $\mathbb{C}^{n_1\times\cdots\times n_d}$ indexed by the tuple $(\alpha_1,\ldots,\alpha_d)\in[n_1]\times\cdots\times[n_d]$.     
\end{definition}

\noindent Intuitively, the MPS parametrization gives an efficient handle on a physically relevant slice of Hilbert space, defined by only $(n_1+\cdots +n_d)r^2$ parameters rather than the full $\prod_k n_k$ needed to parametrize $\mathbb{C}^{n_1}\otimes\cdots \otimes\mathbb{C}^{n_d}$. The bond dimension then gives a knob for specifying the amount of entanglement in the quantum system, with higher bond dimension corresponding to higher entanglement.

An important question in quantum learning theory is \emph{MPS tomography}: given the ability to \emph{measure} multiple copies of an MPS $\ket{\psi}$, output a description of an MPS $\ket{\psi'}$ for which $\ket{\psi'}$ and $\ket{\psi}$ are close, for instance, in the sense of quantum fidelity $|\braket{\psi'|\psi}|^2$. 

This question is closely related to the problem of tensor ring decomposition, but the key difference lies in the access model: whereas in TR decomposition we assume access to the entries of the tensor, which translates here to entries of the quantum state $\ket{\psi}$, in MPS tomography we are only allowed to perform \emph{quantum measurements} of $\ket{\psi}$. For simplicity, here we focus on \emph{projective} measurements:

\begin{definition}[Projective measurements]
    A \emph{projection-valued measure (PVM)} is specified by a finite collection of orthogonal projectors $\{P_1,\dots,P_m\}$ satisfying $P_j P_k =  P_k 1_{\{j=k\}}$ and $\sum_{k=1}^m P_k = \bfI$. Given a state $\ket{\psi}$, measuring $\ket{\psi}$ with this PVM produces an outcome $k \in [m]$ with probability $\Pr(k) = \langle \psi | P_k |\psi\rangle$. After obtaining outcome $k$, the post-measurement state becomes $P_k \ket{\psi}/\sqrt{\langle\psi| P_k |\psi\rangle}$.
    If $\ket{\psi}\in\mathbb{C}^{n_1\times\cdots\times n_d}$, this measurement is said to be \emph{local} if each $P_k$ decomposes as a tensor product of projectors acting independently on $\mathbb{C}^{n_1}, \ldots,\mathbb{C}^{n_d}$, such that for all but a constant number of modes, the projectors all act as the identity on that mode.
\end{definition}

\noindent As quantum measurement results in \emph{probabilistic} access to the tensor given by $\ket{\psi}$, one must perform measurements on multiple copies of the unknown state to estimate it with sufficient accuracy. Beyond this probabilistic aspect, how does the quantum measurement model compare to the access model studied in this paper? The latter corresponds to the model where one gets direct query access to entries of the MPS $\ket{\psi}$. In contrast, with quantum measurements, performing enough measurements to obtain a sufficiently accurate estimate of any prescribed entry of $\ket{\psi}$ is prohibitively expensive because the average entry has squared magnitude $1/\prod_k n_k$ which is exponentially small in $d$.

Instead of getting access to individual entries of the tensor, one kind of information about $\ket{\psi}$ that can be obtained from a small number of measurements is \emph{few-body marginals} of $\ket{\psi}$.
\begin{definition}[Partial trace]
    Given a state $\ket{\psi}\in\mathbb{C}^{n_1\times\cdots\times n_d}$ and a subset $S\subseteq[d]$, the \emph{partial trace} $\mathrm{tr}_{S^c}(\ketbra{\psi})$ is the operator acting on $\otimes_{k\in S} \mathbb{C}^{n_k}$ with $(\balpha,\bbeta)$-th entry given by 
    \begin{equation}
        \sum_{\bgamma\in \prod_{k\not\in S}[n_k]} \psi_{\balpha\bgamma} \psi^*_{\bbeta\bgamma}\,,
    \end{equation}
    for every $\balpha,\bbeta\in\prod_{k\in S}[n_k]$, where $\balpha\bgamma$ is shorthand for the tuple in $\prod^d_{k=1}[n_k]$ whose $k$-th entry is $\alpha_k$ if $k\in S$ and $\gamma_k$ otherwise, and $\bbeta\bgamma$ is defined analogously. This partial trace is the \emph{marginal of $\ket{\psi}$ on subsystem $S$}.

    If $\ket{\psi}$ is an MPS defined by cores $\mbcQ_1,\ldots,\mbcQ_d$, then the $(\balpha,\bbeta)$-th entry is 
    \begin{equation}
        \sum_{\bgamma\in\prod_{k\not\in S}[n_k]} \mathrm{tr}\!\left\{\bfQ^{((\alpha\gamma)_1)}_1\bfQ^{((\alpha\gamma)_2)}_2\cdots \bfQ^{((\alpha\gamma)_d)}_d\right\} \cdot \mathrm{tr}\!\left\{\bfQ^{((\beta\gamma)_1)}_1\bfQ^{((\beta\gamma)_2)}_2\cdots \bfQ^{((\beta\gamma)_d)}_d\right\}^*\,. \label{eq:mpsrdmentry}
    \end{equation}
\end{definition}

\begin{remark}[Phase ambiguity]
    The recovery of an MPS is subject to an additional phase ambiguity beyond the gauge invariance in~\eqref{eq:gauge-invariance}. Specifically, since quantum states $|\psi\rangle$ and $\lambda\,|\psi\rangle$ are physically equivalent for any $\lambda\in \bbS^1$, the MPS parameters can only be recovered up to a global \emph{phase} factor $\lambda$. This is inherent to the quantum setting and does not affect the practical utility of the recovery, since all physically meaningful quantities (e.g., expectation values and quantum fidelity) are invariant under global phase.
\end{remark}

As the partial trace $\mathrm{tr}_S(\ketbra{\psi})$ is a matrix in dimension $\prod_{k\in S}[n_k]$, if $|S| = O(1)$ then its entries have average squared magnitude at least $1/\mathrm{poly}(\max_k n_k)$ and can thus be estimated to sufficient accuracy using only polynomially many measurements. Furthermore, it is a standard fact~\cite[Section 8.4.2]{nielsen2010quantum} that such marginals can be estimated using only local, projective measurements.

Next, we explain how  \blostr can be applied to reconstruct the parameters of a generic matrix product state based on $O(d)$ many \emph{$3$-body marginals} for $\ket{\psi}$. 

\begin{lemma}
    \label{lemma:mps}
    Suppose the bond dimension $r \le \min_{k\in[d]} \sqrt{n_k}$. For any consecutive $j \prec k\in[d]$, $h\in[d]\backslash\{j,k\}$ and unit vector $v\in\mathbb{C}^{n_h}$, given access to the $3$-body marginal $\brho^S \triangleq \mathrm{tr}_{S^c}(\ketbra{\psi})$ for $S = \{h,j,k\}$, one can construct the matrix
    \begin{equation}
        \bfM^{h,v}_{j,k} \triangleq \bfQ_{k\left<1\right>} \bigl(\bfI_r \otimes \bfR^{h,v}_{j,k}\bigr) \bfQ^*_{j[1]}\,\qquad \text{for} \qquad \bfR^{k,v}_{j,k} \triangleq \sum^{n_h}_{\alpha_h , \alpha'_h = 1} v_{\alpha_h} v^*_{\alpha'_h} \sum_{\balpha_{\backslash h}} \tau_{\balpha'}^* \bfR^{j,k;\balpha}\,, \label{eq:Mkv}
    \end{equation}
    where $\balpha_{\backslash h}$ ranges over $\prod_{\ell\in\{k+1,\ldots,j-1\}\backslash \{h\}} [n_\ell]$, $\balpha \in \prod_{k+1\le \ell\le j-1}[n_\ell]$ (resp. $\balpha'$) denotes the concatenation of $\alpha_h$ (resp. $\alpha'_h$) and $\balpha_{\backslash h}$ defined in the natural way, $\tau_{\balpha'}$ is the entry of $\mbcT$ which is indexed by $1$ in the $j$-th and $k$-th modes, and $\balpha'$ in the remaining modes, and $\bfR^{j,k;\balpha} = \prod_{\ell = k+1}^{j-1} \bfQ^{(\alpha_\ell)}_\ell$. 
\end{lemma}

\begin{remark}[On the dimension assumption in the MPS setting]\label{remark:dimension-assumption-mps}
The condition $n_k \ge r^2$ appearing in our MPS recovery result is equivalent to $r \le \sqrt{n_k}$ for all $k$. We view this as a moderate-bond-dimension regime. While many quantum systems have small single-site physical dimension, a standard operation in tensor-network methods is to block a constant number of neighboring sites into one effective site, which replaces the local dimension $q$ by $q^m$. Hence the assumption can often be satisfied after coarse-graining, without essentially changing the underlying periodic-MPS structure.
\end{remark}

Indeed, one can obtain $\bfM_{j,k}^{h,v}$ by first contracting  $\brho^S$ along the $h$-th mode to form an $n_jn_k\times n_jn_k$ matrix, then reshaping it so that the entry indexed by $((\alpha_k,\alpha_j),(\alpha_k',\alpha_j'))$ is moved to $((\alpha_k,\alpha_k'),(\alpha_j,\alpha_j'))$. Finally, we take the entries $((\alpha_k,1),(\alpha_j,1))$ for all $\alpha_j\in[n_j]$ and $\alpha_k\in[n_k]$.
A detailed procedure is provided in Appendix~\ref{sec:proof-lemma}. 

In the definition of $\bfM^{h,v}_{j,k}$ in Eq.~\eqref{eq:Mkv}, the operator $\bfR^{h,v}_{j,k}$ is invertible with probability 1. We can thus proceed as in the proof of Theorem~\ref{thm:order-d}, but with $\bfR^{\balpha}, \bfR^{\bbeta}, \bfR^{\balpha'}, \bfR^{\bbeta'}$ therein replaced with $\bfR^{h,v}_{j,k}$ for four different choices of $(h, v)$, and with $\bfT(:,\balpha_{\rm mid},:)$ therein replaced by $\bfM^{h,v}_{j,k}$. The argument in Steps 1 to 3 of Section~\ref{sec:procedure} proceeds entirely analogously, including the computation of $\hat{\bfQ}_{k\left<1\right>} = \bfE(\Pi(\hat{\bfK}))^{-1}$ in Step 3, with $\bfE$ (resp. $\bfE'$) taken to be an eigenbasis for the nonzero spectrum of $\bfM^{h_1,v_1}_{j,k} (\bfM^{h_2,v_2}_{j,k})^{\dagger}$ (resp. $\bfM^{h_3,v_3}_{j,k} (\bfM^{h_4,v_4}_{j,k})^{\dagger}$), and $\bfF, \hat{\bfK}$ constructed as in Step 2. We take eigendecompositions $\bfR^{h_1,v_1}_{j,k}(\bfR^{h_2,v_2}_{j,k})^{-1} = \bfU \boldsymbol{\Lambda} \bfU^{-1}$ and $\bfR^{h_3,v_3}_{j,k} (\bfR^{h_4,v_4}_{j,k})^{-1} = \bfV \boldsymbol{\Lambda}'\bfV^{-1}$. By a straightforward adaptation of Lemma~\ref{lemma:eigenvalues}, $\bfE = \bfQ_{1\left<1\right>} (\bfI_r\otimes \bfU) \Pi(\bfK)$ for some block-diagonal $\bfK$, and similarly for $\bfE'$. This means that $\bfF \triangleq \bfE^\dagger \bfE'$ also satisfies Eq.~\eqref{eq:F}. We then apply Eq.~\eqref{eq:hatq-1} to obtain $\hat{\mbcQ}_k$, and set $\tilde{\bfQ}_{j[1]}=\bigl(\hat{\bfQ}_{k\left<1\right>}^\dagger \bfM_{j,k}^{h,v}\bigr)^*$.

Following a calculation similar to the one in Appendix~\ref{sec:proof}, one can show that the obtained pair $\{\tilde{\mbcQ}_j,\hat{\mbcQ}_k\}$ satisfies
\begin{align}
    \label{eq:ij-pairs}
    \tilde{\mbcQ}_j = \mbcQ_j \times_2 \bfZ_k^{-1} \times_3 \bfX_k^*,\qquad 
    \hat{\mbcQ}_k = \mbcQ_k \times_2 \bfX_k^{-1}\times_3 \bfY_k^*,\qquad 
    \text{for some}\ \bfX_k,\bfY_k,\bfZ_k\in\GL(r,\bbC).
\end{align}

Therefore, we can get $d$ pairs of cores in the form of Eq.~\eqref{eq:ij-pairs}  by performing $(d-1)$ permutations over $\mbcT$. We then match the gauge of those sharing the same subscript to obtain a set of compatible cores. In particular, notice that 
\begin{align*}
    \bigl\{\tilde{\mbcQ}_1 \times_2 (\bfX_1^{-1}\bfZ_2),\ \tilde{\mbcQ}_2 \times_2 (\bfX_2^{-1}\bfZ_3),\ \ldots,\ \tilde{\mbcQ}_d \times_2 (\bfX_d^{-1}\bfZ_1)\bigr\}
\end{align*}
forms a set of cores; and for $k\in[d]$,
\begin{align*}
    \tilde{\bfQ}_{k[1]}^\dagger \hat{\bfQ}_{k[1]} &= \left(\bfQ_{k[1]} (\bfX_{k+1}^* \otimes \bfZ_{k+1}^{-1})^*\right)^\dagger \left(\bfQ_{k[1]}(\bfY_k^* \otimes \bfX_1^{-1})^*\right)\\
    &= (\bfX_{k+1}^{-1} \otimes \bfZ_{k+1}^{*}) (\bfY_k \otimes \bfX_k^{*-1}) \\
    &= (\bfX_{k+1}^{-1} \bfY_k \otimes (\bfX_k^{-1} \bfZ_{k+1})^*),
\end{align*}
which implies that, up to some scalar, $(\bfX_k^{-1}\bfZ_{k+1})^{*}$ can be identified with the top-left $r\times r$ sub-matrix of $\tilde{\bfQ}_{k[1]}^\dagger \hat{\bfQ}_{k[1]}$, which we denote by $\bfN_k$. We then rescale $\hat{\mbcQ}_k\times_2 \bfN_k$ to ensure consistent scalar propagation. We summarize the preceding discussion in the following theorem.

\begin{theorem}[\blostr recovers MPS structure]\label{thm:mps}
    Suppose $\ket{\psi}\in\bbC^{n_1} \otimes\cdots\otimes \bbC^{n_d}$ is an MPS with bond dimension $r$ in the sense of \eqref{eq:MPS}, where $d\geq 3$, $r\geq 2$, and $n_k \geq r^2$ for all $k\in[d]$. Assume that the elements of $\mbcQ_k\in\mathbb{C}^{n_k\times r\times r}$ are randomly drawn from a measure $\mu$ that is absolutely continuous with respect to the standard Euclidean measure, and suppose we are given access to the $3$-body marginals of $\ket{\psi}$. Then with probability one, the parameters $\mbcQ_1, \mbcQ_2,\dots, \mbcQ_d$ of the MPS can be recovered by Algorithm~\ref{alg:mps} up to gauge invariance in \eqref{eq:gauge-invariance} and a global phase factor $\lambda\in\bbS^1$.
\end{theorem}

\begin{algorithm}[ht]
    \caption{\blostr for MPS recovery}
    \label{alg:mps}
    \KwIn{$\bigl\{\bfM_{j,k}^{h_\iota,v_\iota}:j\in[d],\ k=j+1;\ \iota=1,2,3,4\bigr\}$, measurements of $\ket{\psi}$}
    \KwOut{$\check{\mbcQ}_k\in\bbC^{n_k\times r\times r}$, $k\in[d]$, solution to \eqref{eq:MPS}}

    \For{$j=1$ \KwTo $d$}{
        Perform lines 2--6 in Algorithm~\ref{alg:order-d} with $\bfT(:,\balpha,:),  \bfT(:,\bbeta,:), \bfT(:,\balpha',:), \bfT(:,\bbeta',:)$ replaced by $\{\bfM_{j,k}^{h_\iota,v_\iota}\}_{\iota=1}^4$ to obtain $\hat{\mbcQ}_k$\;
        Obtain $\tilde{\mbcQ}_j$ by computing $\tilde{\bfQ}_{j[1]}=(\hat{\bfQ}_{k\langle1\rangle}^\dagger\bfM_{j,k}^{h_1,v_1})^*$ and reshaping\;
    }
    
    \For{$k=1$ \KwTo $d$}{
        Let $\bfN_k$ be the top-left $r\times r$ sub-matrix of $\tilde{\bfQ}_{k[1]}^\dagger \hat{\bfQ}_{k[1]}$\;
        $\check{\mbcQ}_k \gets \tilde{\mbcQ}_k \times_2 \bfN_k^{*}$\;
    }
    Set $\check{\boldsymbol{\mathcal{Q}}}_1  \gets \check{\boldsymbol{\mathcal{Q}}}_1 / \|\ket{\psi'}\|$, where $\ket{\psi'} = \sum_{\alpha_1,\ldots,\alpha_d} \mathrm{tr}\!\left\{ \check{\mathbf{Q}}^{(\alpha_1)}_1 \check{\mathbf{Q}}^{(\alpha_2)}_2  \cdots \check{\mathbf{Q}}^{(\alpha_d)}_d\right\} \ket{\alpha_1,\ldots,\alpha_d}$\;
\end{algorithm}

\begin{remark}
    We note that in the usual model of quantum state learning, one does not get exact access to the $3$-body marginals of $\ket{\psi}$, but instead \emph{approximate} access from performing measurements on finitely many copies of the state. This ultimately incurs some stochastic errors which must be accounted for. Given the robustness of \blostr, this should not affect the performance of our algorithm, but we defer a quantitative analysis of this error propagation and a careful quantification of the sample complexity of our MPS tomography protocol to future work.
\end{remark}

\begin{remark}
    There have been numerous prior works on MPS tomography which target slightly different guarantees. Efficient algorithms were given in~\cite{landon2010efficient,cramer2010efficient} for the problem of \emph{learning a quantum circuit for preparing the MPS} given the ability to perform \emph{global} (but efficient) measurements of copies of the state. The work of~\cite{cramer2010efficient} also proposed a heuristic method for learning the circuit using only local measurements, but to our knowledge no formal such guarantee is known in the literature. Indeed, for \emph{worst-case} matrix product states, it is in general not possible to recover them from local measurements alone (e.g., CAT states are indistinguishable from the maximally mixed state unless one considers global marginals), and our recovery guarantee above sidesteps this by considering \emph{generic} MPS's that are drawn from an absolutely continuous measure. More generally one might hope for learning algorithms for \emph{injective} MPS's, for which the works of~\cite{landau2015polynomial,arad2017rigorous} give algorithms for reconstructing an approximation to the MPS, with the caveat that the output can have \emph{larger} bond dimension than the original state. We defer a detailed discussion of this literature to~\cite{anshu2024survey}.
\end{remark}

\subsection{Learning Polynomial Transformations}
\label{sec:main_moment}

Beyond the MPS interpretation from physics, the TR decomposition also admits a statistical viewpoint through its connection to \emph{moment tensors}. Moment tensors are fundamental for capturing higher-order correlations in latent-variable models, and their decomposition has long been a central tool in statistical learning. Building upon an observation from~\cite{chen2023learning}, we show that certain polynomial transformations of Gaussian distributions can be learned from their moments up to degree $d$ by solving a (noisy) TR decomposition problem.

\begin{definition}[Polynomial transformations]\label{def:polytransform}
    A \emph{cyclic quadratic transformation} is a distribution $q$ over $\mathbb{R}^{n_1+\cdots+n_d}$ specified by cores $\mbcQ_1 \in \mathbb{R}^{n_1\times r\times r},\ldots,\mbcQ_d \in \mathbb{R}^{n_d\times r\times r}$ through the following pushforward construction. To sample from $q$, sample $x_1,\ldots,x_d\sim\mathcal{N}(0,\bfI_r)$ and output the vector $z$ consisting of $d$ blocks of coordinates of size $n_1,\ldots,n_d$ respectively, such that the $\alpha_k$-th entry in the $k$-th block is given by
    \begin{equation}
        \label{eq:cyclic-quadratic}
        z^{(\alpha_k)}_k \triangleq \bfx^*_k \bfQ^{(\alpha_k)}_k \bfx_{k+1}\,,
    \end{equation}
    where for $k = d$ we define $\bfx_{d+1} = \bfx_1$.
    
    The problem of \emph{parameter learning} for such distributions is the following: given $N$ i.i.d. samples drawn from $q$, output estimates $\hat{\mbcQ}_k$ for the cores $\mbcQ_k$ up to gauge invariance.
    
\end{definition}

\noindent In this polynomial transformation problem, the condition \(n_k \ge r^2\) is natural in the regime where the observed output dimension is large relative to the latent dimension \(r\). A cyclic quadratic transformation can be viewed as a stylized generative model that maps a latent space of dimension \(rd\) to an observed space of dimension \(n_1+\cdots+n_d\), with \(n_1,\ldots,n_d\) naturally playing the role of output dimensions. The output space is divided into $d$ patches, each computed by some degree-2 polynomials in the latents, with some consistency between any $(k-1)$-st and $k$-th patch enforced by the fact that they share dependence on the $x_k$ block of latents.

The following lemma, whose proof is immediate from Wick calculus, establishes the connection between parameter learning of cyclic quadratic transformations and TR decomposition.

\begin{lemma}
    Let $q$ be a polynomial transformation specified by $\mbcQ_1,\ldots,\mbcQ_d$ as in Definition~\ref{def:polytransform}. 
    Then for any $\alpha_1\in [n_1],\ldots,\alpha_d\in[n_d]$, for $z\sim q$,
    \begin{equation}
        \mathbb{E}[z^{(\alpha_1)}_1\cdots z^{(\alpha_d)}_d] = \mathrm{tr}\!\left\{\bfQ_1^{(\alpha_1)}\bfQ_2^{(\alpha_2)}\cdots\bfQ_d^{(\alpha_d)}\right\}\,.
    \end{equation}
\end{lemma}

\noindent In other words, if one looks at the order-$d$ \emph{moment tensor} whose entries correspond to all degree-$d$ moments of $q$, the entries of this tensor contain the tensor $\mbcT$ whose TR decomposition is given by the cores $\mbcQ_1,\ldots,\mbcQ_d$ defining $q$. As a consequence, this simple family of generative models can be parameter learned using the method of moments, in conjunction with \blostr. Furthermore, as noted in our discussion at the end of Section~\ref{sec:procedure}, \blostr only needs to inspect $O(r^2 \sum_j n_j)$ total entries, as opposed to all $\prod_k n_k$ entries of $\mbcT$, so once the relevant moments corresponding to those entries have been estimated to sufficient accuracy, the runtime of our learning algorithm scales only polynomially in the data dimension, rather than exponentially in the number $d$ of patches.

\begin{figure}[ht]
    \centering
    \includegraphics[width=\linewidth]{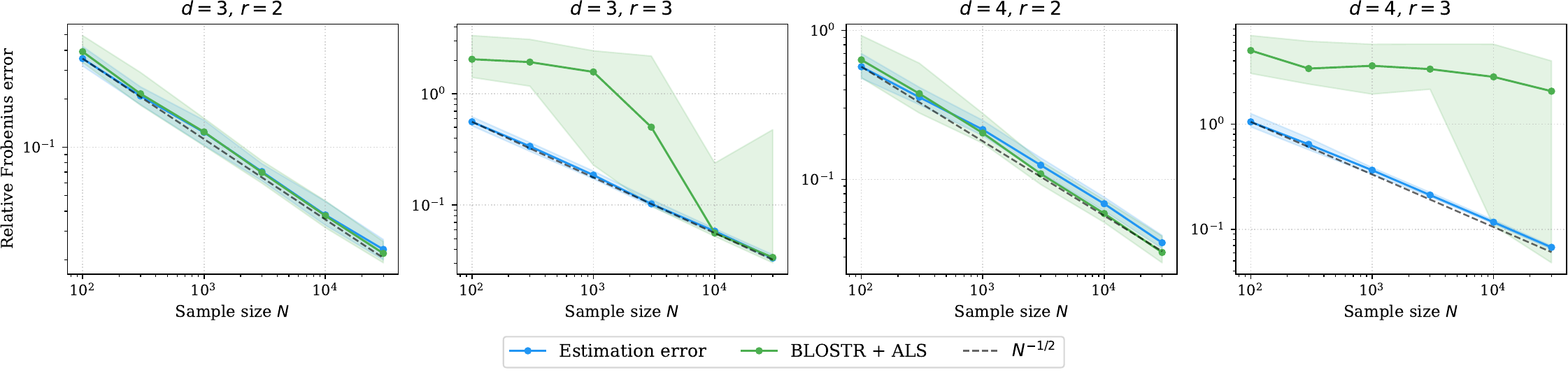}
    \caption{Relative reconstruction error versus sample size $N$ for the polynomial transformation model~\eqref{eq:cyclic-quadratic}, with $d\in\{3,4\}$ and $r\in\{2,3\}$. The cores $\{\mbcQ_k\}_{k\in[d]}$ have i.i.d.\ $\mathcal{N}(0,1)$ entries. For each trial, $N$ i.i.d.\ samples are drawn and the population moment tensor $\mathcal{T}$ is estimated by sample averaging to obtain $\hat{\mathcal{T}}$. Blue: estimation error $\|\mathcal{T}-\hat{\mathcal{T}}\|_F/\|\mathcal{T}\|_F$.  Green: reconstruction error after \blostr initialization followed by ALS refinement. Shaded bands indicate the interquartile range over 100 independent trials; the dashed line shows the $N^{-1/2}$ reference slope. }
    \label{fig:polytransform}
\end{figure}

We provide a numerical illustration of the method-of-moments approach described above. We generate synthetic data from the model in Definition~\ref{def:polytransform} with $d \in \{3, 4\}$, $n_k = n = 12$ for all $k \in [d]$, and TR rank $r \in \{2, 3\}$. The cores $\{\mbcQ_k\}_{k \in [d]}$ have entries drawn i.i.d.\ from $\mathcal{N}(0, \sigma_s^2)$ with $\sigma_s = 1$. For each trial, we generate $N$ i.i.d.\ samples from the transformation. Specifically, for the $j$-th sample, we draw $\bfx_1^{j}, \ldots, \bfx_d^{j} \sim \mathcal{N}(0, I_r)$ independently and compute the output $z^{(j)}$ via~\eqref{eq:cyclic-quadratic}, where $z_k^{(\alpha_k, j)} = (\bfx_k^j)^\ast \bfQ_k^{(\alpha_k)} \bfx_{k+1}^j$ denotes the 
$\alpha_k$-th entry of the $k$-th block of the $j$-th sample. We then estimate the moment tensor $\boldsymbol{\mathcal{T}}$ via sample averaging, i.e.,
\begin{equation}
    \hat{T}(\alpha_1, \ldots, \alpha_d) = \frac{1}{N} \sum_{j=1}^N 
    z_1^{(\alpha_1, j)} \cdots z_d^{(\alpha_d, j)},
\end{equation}
and apply \blostr to $\hat{\boldsymbol{\mathcal{T}}}$ to recover the cores. We report the relative reconstruction error  averaged over 100 independent trials as a function of the sample size $N$. The results are shown in Figure~\ref{fig:polytransform}. Our algorithm recovers the TR parameters at the optimal $N^{-1/2}$ rate when the TR rank $r$ is small; for larger $r$, convergence is slower but the trend is consistent.

\begin{remark}\label{remark:symmetric_poly_transform}    In~\cite{chen2023learning}, the authors provided an algorithm for learning non-cyclic quadratic transformations in which every coordinate is a generic quadratic form in \emph{all} of the seed variables. They observed that the third-order cumulants of such distributions are precisely the entries of a symmetric tensor with bounded TR rank, as in the setting of our Section~\ref{sec:symmetric}, but with the additional twist that the matrices $\mbcQ^{(\alpha)}$ in the symmetric TR core $\mbcQ$ are \emph{symmetric} matrices. They provided a highly involved estimator based on the sum-of-squares algorithm for this version of the TR decomposition, and we leave it as an interesting open question to obtain a simple, \blostr-style algorithm in this setting.
\end{remark}

\section{Discussion}\label{sec:discussion}

This work provides the first finite-step algorithm for tensor ring (TR) decomposition. By introducing a blockwise simultaneous diagonalization method, we establish a deterministic procedure \blostr that recovers TR cores exactly from sparse observations, thereby uncovering the algebraic structure of the TR format. This result places TR decomposition on par with classical tensor models such as Tucker \cite{de2000best,de2000multilinear,zhang2018tensor} and tensor-train \cite{oseledets2011tensor,zhou2022optimal}, for which finite-step methods have long been known.

Beyond this theoretical advance, we extend the method to symmetric TR decomposition, showing how symmetry can be leveraged to reduce parameter complexity. We further address noisy and imperfect observations through a robust recovery scheme that combines \blostr with alternating least squares refinement, demonstrating empirically that principled initialization dramatically accelerates convergence and improves stability compared to standard random-initialization methods.

We acknowledge several limitations of the proposed method. First, the algorithm assumes the TR rank $r$ is known a priori; when it is not, rank estimation is required as a preprocessing step. Second, the identifiability guarantee relies on the genericity assumption on the cores, which may not hold for structured or sparse cores arising in certain applications. Third, while the dimension condition \(n_k \ge r^2\) is satisfied in many practical settings and, in quantum/tensor-network applications, can often be enforced by blocking, it may be restrictive in regimes where the bond dimension $r$ is large relative to the mode sizes $n_k$. We leave the development of methods for these more challenging regimes as future work.

Looking ahead, several broader research directions emerge. A key challenge is to generalize beyond the equal-rank setting, where tensor cores may have heterogeneous dimensions. Another is to develop adaptive procedures that do not require prior knowledge of the TR rank, a critical step toward making these methods broadly applicable in practice. More ambitiously, extending finite-step guarantees to noisy settings, with explicit error bounds and robustness analyses, would bridge the gap between algebraic exactness and statistical efficiency.  

Finally, our blockwise simultaneous diagonalization-based framework resonates with a growing body of work at the interface of machine learning, physics, and applied mathematics. The connections between TR decomposition, tensor network models for many-body quantum systems such as matrix product states, and moment-based methods for learning latent variable distributions suggest fertile ground for interdisciplinary progress. We believe that the tools introduced here will inspire new algorithmic paradigms for structured representation learning and the principled study of high-dimensional high-order data.

\paragraph{Code availability and reproducibility}
Code for implementing \blostr and reproducing the results in this paper is available at \href{https://github.com/chenh77/blostr.git}{\tt github.com/chenh77/blostr}.

\subsection*{Acknowledgments}

The authors thank Anurag Anshu, Luke Coffman, and Quynh Nguyen for insightful discussions about matrix product state tomography. The authors also thank the anonymous referees for their careful reading and helpful comments.

\bibliography{reference-new}


\clearpage

\appendix
\appendixpage
\setcounter{lemma}{0}
\setcounter{table}{0}
\setcounter{figure}{0}
\setcounter{equation}{0}
\setcounter{algocf}{0}
\renewcommand{\thelemma}{S\arabic{lemma}}
\renewcommand{\thetable}{\Alph{section}\arabic{table}}
\renewcommand{\thefigure}{\Alph{section}\arabic{figure}}
\renewcommand{\theequation}{\Alph{section}.\arabic{equation}}
\renewcommand{\thealgocf}{\Alph{section}\arabic{algocf}}

\section{Proofs of Theorems}\label{sec:proof}

\subsection{Auxiliary Lemmas}

The first lemma is a classic conclusion on the Moore-Penrose inverse \cite{ben2003generalized}. 
\begin{lemma}\label{lemma:mp-inverse}
    Suppose $\bfA$ and $\bfB$ are two matrices with compatible dimensions. If $\bfA$ has full column rank, and $\bfB$ has full row rank, then $(\bfA\bfB)^\dagger=\bfB^\dagger\bfA^\dagger$. 
\end{lemma}

The following two lemmas explore some useful properties of the operation $\Pi_{r_1,r_2}$ defined in Section~\ref{sec:preliminary}. We use $\circ$ for the function composition.

\begin{lemma}\label{lemma:kronecker}
    For $r_1,r_2\geq2$, let $\Pi_{r_1,r_2}$ be as defined in \eqref{eq:pi}. The following holds:
    \begin{enumerate}[label=(\roman*)]
        \item $\Pi_{r_1,r_2}\circ\Pi_{r_2,r_1}(\bfX)=\Pi_{r_2,r_1}\circ\Pi_{r_1,r_2}(\bfX)=\bfX$;
        \item Assume that $\bfX\in\bbC^{r_1\times r_1}$ and $\bfY\in\bbC^{r_2\times r_2}$, then $\Pi_{r_1,r_2}(\bfX\otimes\bfY)=\bfY\otimes\bfX$.
    \end{enumerate} 
\end{lemma}

\begin{proof}[Proof of Lemma~\ref{lemma:kronecker}]~
\begin{enumerate}[label=(\roman*)]
    \item For $j_1,k_1\in[r_1]$ and $j_2,k_2\in[r_2]$, we have 
    \begin{align*}
        \left(\Pi_{r_1,r_2}\circ\Pi_{r_2,r_1}(\bfX)\right)_{(j_2-1)r_1+j_1,(k_2-1)r_1+k_1}&=\left(\Pi_{r_2,r_1}(\bfX)\right)_{(j_1-1)r_2+j_2,(k_1-1)r_2+k_2}\\
        &=\bfX_{(j_2-1)r_1+j_1,(k_2-1)r_1+k_1},\\
        \left(\Pi_{r_2,r_1}\circ\Pi_{r_1,r_2}(\bfX)\right)_{(j_1-1)r_2+j_2,(k_1-1)r_2+k_2}&=\left(\Pi_{r_1,r_2}(\bfX)\right)_{(j_2-1)r_1+j_1,(k_2-1)r_1+k_1}\\
        &=\bfX_{(j_1-1)r_2+j_2,(k_1-1)r_2+k_2}.
    \end{align*}
    Therefore, $\Pi_{r_1,r_2}\circ\Pi_{r_2,r_1}(\bfX)=\Pi_{r_2,r_1}\circ\Pi_{r_1,r_2}(\bfX)=\bfX$.
    \item For $j_1,k_1\in[r_1]$ and $j_2,k_2\in[r_2]$, we have
    \begin{align*}
        \left(\Pi_{r_1,r_2}(\bfX\otimes\bfY)\right)_{(j_2-1)r_1+j_1,(k_2-1)r_1+k_1}&=(\bfX\otimes\bfY)_{(j_1-1)r_2+j_2,(k_1-1)r_2+k_2}\\
        &=\bfX_{j_1,k_1}\bfY_{j_2,k_2}\\
        &=\bfY_{j_2,k_2}\bfX_{j_1,k_1}\\
        &=(\bfY\otimes\bfX)_{(j_2-1)r_1+j_1,(k_2-1)r_1+k_1}.
    \end{align*}
    Therefore, $\Pi_{r_1,r_2}(\bfX\otimes\bfY)=\bfY\otimes\bfX$.
\end{enumerate}
\end{proof}

\begin{lemma}\label{lemma:pi}
    Assume that $\bfK,\bfH\in\bbC^{r_1r_2\times r_1r_2}$ are blockwise diagonal matrices, $\bfK=\text{diag}(\bfK_1,\dots,\bfK_{r_1})$, $\bfH=\text{diag}(\bfH_1,\dots,\bfH_{r_1})$, with $\bfK_t\in\bbC^{r_2\times r_2}$, $\bfH_t\in\bbC^{r_2\times r_2}$ for $t\in[r_1]$. Then $\Pi_{r_1,r_2}(\bfK)\,\Pi_{r_1,r_2}(\bfH)=\Pi_{r_1r_2}(\bfK\bfH)$. In particular, $\Pi_{r_1,r_2}(\bfK^{-1})=(\Pi_{r_1,r_2}(\bfK))^{-1}$.
\end{lemma}

\begin{proof}[Proof of Lemma~\ref{lemma:pi}]
    For $\ell\in[r_1]$, $j,k\in[r_2]$, notice that
    \begin{align*}
        & \left(\Pi_{r_1,r_2}(\bfK)\Pi_{r_1,r_2}(\bfH)\right)_{(j-1)r_1+\ell,(k-1)r_1+\ell}\\
        =&\sum_{m=1}^{r_2}\Pi_{r_1,r_2}(\bfK)_{(j-1)r_1+\ell,(m-1)r_1+\ell}\cdot\Pi_{r_1,r_2}(\bfH)_{(m-1)r_1+\ell,(k-1)r_1+\ell}\\
        =&\sum_{m=1}^{r_2}(\bfK_\ell)_{j,m}(\bfH_\ell)_{m,k}\\
        =&((\bfK\bfH)_\ell)_{j,k}\\
        =&(\Pi_{r_1,r_2}(\bfK\bfH))_{(j-1)r_1+\ell,(k-1)r_1+\ell},
    \end{align*}
    and that 
    \begin{align*}
        \left(\Pi_{r_1,r_2}(\bfK)\Pi_{r_1,r_2}(\bfH)\right)_{(j_2-1)r_1+j_1,(k_2-1)r_1+k_1}=0,\ \text{for}\ j_1\neq k_1.
    \end{align*} 
    Therefore $\Pi_{r_1,r_2}(\bfK)\Pi_{r_1,r_2}(\bfH)=\Pi_{r_1,r_2}(\bfK\bfH)$. Then we have 
    \begin{align*}
        \Pi_{r_1,r_2}(\bfK)\Pi_{r_1,r_2}(\bfK^{-1})=\Pi_{r_1,r_2}(\bfK\bfK^{-1})=\bfI_{r_1r_2}=\Pi_{r_1,r_2}(\bfK^{-1}\bfK)=\Pi_{r_1,r_2}(\bfK^{-1})\Pi_{r_1,r_2}(\bfK). 
    \end{align*}
    That is, $(\Pi_{r_1,r_2}(\bfK))^{-1}=\Pi_{r_1,r_2}(\bfK^{-1})$. 
\end{proof}

\subsection{Proof of Theorem~\ref{thm:order-d}}\label{sec:proof-theorem1}

We now provide the detailed proof, following the outline given in Section~\ref{sec:procedure}. Throughout, we work on the probability--one event guaranteed by the absolute continuity of $\mu$, on which all random matrices appearing below are invertible and have simple spectra. In particular, every $\bfR^{\balpha}$ is invertible, and $\bfR^{\balpha}(\bfR^{\bbeta})^{-1}$ and $\bfR^{\balpha'}(\bfR^{\bbeta'})^{-1}$ are diagonalizable with distinct eigenvalues. We fix an ordering of eigenvectors by ordering eigenvalues, so that the two spectral probes are aligned up to blockwise scaling; remaining ambiguities are absorbed into the gauge matrices below.

\paragraph{Step 1 (two spectral probes)}
Since $\bfR^{\balpha}\in\bbC^{r\times r}$ for all $\balpha\in[n_2]\times\cdots\times[n_{d-1}]$, and the elements drawn from $\mu$ that is absolutely continuous with respect to the standard Euclidean measure, $\bfR^{\balpha}$ is invertible with probability one. Let $\bfR^{\balpha}(\bfR^{\bbeta})^{-1}=\bfU\boldsymbol{\Lambda}\bfU^{-1}$ and $\bfR^{\balpha'}(\bfR^{\bbeta'})^{-1}=\bfV\boldsymbol{\Lambda}'\bfV^{-1}$ be the eigendecomposition of $\bfR^{\balpha}(\bfR^{\bbeta})^{-1}$ and $\bfR^{\balpha'}(\bfR^{\bbeta'})^{-1}$, respectively.  
By Lemma~\ref{lemma:eigenvalues}, the eigenvectors corresponding to the nonzero eigenvalues of $\bfT(:,\balpha,\Gamma_{\balpha})\bfT(:,\bbeta,\Gamma_{\bbeta})^\dagger$ have the form $\bfE=\bfQ_{1\left<1\right>}\left(\bfI_{r}\otimes\bfU\right)\Pi(\bfK)$, and those of $\bfT(:,\balpha',\Gamma_{\balpha'})\{\bfT(:,\bbeta',\Gamma_{\bbeta'})\}^\dagger$ have the form $\bfE'=\bfQ_{1\left<1\right>}\left(\bfI_{r}\otimes\bfV\right)\Pi(\bfK')$ for some $\bfK=\text{diag}(\bfK_1,\dots,\bfK_{r})$ and $\bfK'=\text{diag}(\bfK_1',\dots,\bfK_{r}')$, where $\bfK_j,\bfK_j'\in\GL(r,\bbC)$, $j\in[r]$. Since $n_1\geq r^2$, $\bfQ_{1\left<1\right>}$ is of full column rank with probability one. Following Lemma~\ref{lemma:mp-inverse}, with probability one, we have
\begin{align*}
    \bfF&=\bfE^\dagger\bfE'=\left(\bfQ_{1\left<1\right>}\left(\bfI_{r}\otimes\bfU\right)\Pi(\bfK)\right)^\dagger\left(\bfQ_{1\left<1\right>}\left(\bfI_{r}\otimes\bfV\right)\Pi(\bfK')\right)\\
    &=\left(\Pi(\bfK)\right)^{-1}(\bfI_{r}\otimes\bfU^{-1})\bfQ_{1\left<1\right>}^\dagger\bfQ_{1\left<1\right>}\left(\bfI_{r}\otimes\bfV\right)\Pi(\bfK')=\Pi(\bfK^{-1})\left(\bfI_{r}\otimes\left(\bfU^{-1}\bfV\right)\right)\Pi(\bfK'),
\end{align*}
where the last equality follows from Lemma~\ref{lemma:pi} and that $\bfQ_{1\left<1\right>}^\dagger\bfQ_{1\left<1\right>}=\bfI_{r^2}$.
Notice that
\begin{align*}
    &\left(\left(\bfI_{r}\otimes\left(\bfU^{-1}\bfV\right)\right)\Pi(\bfK')\right)_{(j_1-1)r+j_2,(k_1-1)r+k_2}\\
    &=\left(\bfI_{r}\otimes\left(\bfU^{-1}\bfV\right)\right)_{(j_1-1)r+j_2,(j_1-1)r+k_2}\left(\Pi(\bfK')\right)_{(j_1-1)r+k_2,(k_1-1)r+k_2}\\
    &=\left(\bfU^{-1}\bfV\right)_{j_2,k_2}\left(\bfK'_{k_2}\right)_{j_1,k_1}.
\end{align*}
Therefore,
\begin{align*}
    & \bfF_{(j_1-1)r+j_2,(k_1-1)r+k_2}\\
    =&\sum_{\ell=1}^{r}\left(\Pi(\bfK^{-1})\right)_{(j_1-1)r+j_2,(\ell-1)r+j_2}\left(\left(\bfI_{r}\otimes\left(\bfU^{-1}\bfV\right)\right)\Pi(\bfK')\right)_{(\ell-1)r+j_2,(k_1-1)r+k_2}\\
    =&\sum_{\ell=1}^{r}\left(\bfK_{j_2}^{-1}\right)_{j_1,\ell}\left(\bfU^{-1}\bfV\right)_{j_2,k_2}\left(\bfK'_{k_2}\right)_{\ell,k_1}\\
    =&\left(\bfU^{-1}\bfV\right)_{j_2,k_2}\cdot\sum_{\ell=1}^{r}\left(\bfK_{j_2}^{-1}\right)_{j_1,\ell}\left(\bfK'_{k_2}\right)_{\ell,k_1},
\end{align*}
from which we obtain
\begin{align*}
    \bfF_{[0:(r-1)]\cdot r+j_2,[0:(r-1)]\cdot r+k_2}=\left(\bfU^{-1}\bfV\right)_{j_2,k_2}\cdot\bfK_{j_2}^{-1}\bfK'_{k_2}.
\end{align*}

\paragraph{Step 2 (block identification / gauge fixing)}
For $j,k,\ell\in[r]$, let $\bfF^{(j,k)}$, $\hat{\bfK}_\ell$ be as defined in Algorithm~\ref{alg:order-d}. Then we have 
\begin{align*}
    \hat{\bfK}_\ell=\bfF^{(1,1)}(\bfF^{(\ell,1)})^{-1}=(\bfU^{-1}\bfV)_{1,1}\bfK_1^{-1}\bfK_1'(\bfK'_1)^{-1}\bfK_\ell(\bfU^{-1}\bfV)_{\ell,1}^{-1}=\frac{(\bfU^{-1}\bfV)_{1,1}}{(\bfU^{-1}\bfV)_{\ell,1}}\cdot\bfK_1^{-1}\bfK_\ell.
\end{align*}

\paragraph{Step 3 (recover $\mbcQ_1$)}
Let $\bfW=\text{diag}\left((\bfU^{-1}\bfV)_{[1:r],1}\right)\in\bbC^{r\times r}$. Then, 
\begin{equation}
\begin{split}\label{eq:hata}
    \hat{\bfQ}_{1\left<1\right>}&=\bfE(\Pi(\hat{\bfK}))^{-1}=\bfE\cdot\Pi(\hat{\bfK}^{-1})\\
    &=\bfE\cdot\left((\bfU^{-1}\bfV)_{1,1}\right)^{-1}\cdot\Pi\left(\begin{pmatrix}
        \bfK_1^{-1} & & \\
        & \ddots & \\
        & & \bfK_{r}^{-1}
    \end{pmatrix}\begin{pmatrix}
        (\bfU^{-1}\bfV)_{1,1}\bfK_1 & & \\
        & \ddots & \\
        & & (\bfU^{-1}\bfV)_{r,1}\bfK_1
    \end{pmatrix}\right)\\
    &=\left((\bfU^{-1}\bfV)_{1,1}\right)^{-1}\cdot\bfE\cdot\Pi(\bfK^{-1})\Pi\left(\bfW\otimes\bfI_{r}\right)\Pi(\bfI_{r}\otimes\bfK_1)\\
    &=\left((\bfU^{-1}\bfV)_{1,1}\right)^{-1}\cdot\bfQ_{1\left<1\right>}\left(\bfI_{r}\otimes\bfU\right)\Pi(\bfK)\Pi(\bfK^{-1})\left(\bfI_{r}\otimes\bfW\right)(\bfK_1\otimes\bfI_{r})\\
    &=\left((\bfU^{-1}\bfV)_{1,1}\right)^{-1}\cdot\bfQ_{1\left<1\right>}\left(\bfK_1\otimes\bfU\bfW\right),
\end{split}
\end{equation}
where the second equation is from Lemma~\ref{lemma:pi}, and the fourth equation is from Lemma~\ref{lemma:kronecker}. Denote $\bfX\equiv (\bfK_1^*)^{-1}$, and $\bfY\equiv ((\bfU^{-1}\bfV)_{1,1})^{-1}\bfU\bfW$. Then for all $\alpha_1\in[n_1]$,
\begin{align*}
    \hat{\bfQ}_1^{(\alpha_1)}&=\left((\bfU^{-1}\bfV)_{1,1}\right)^{-1}\bfK_1^*\cdot\bfQ_1^{(\alpha_1)}\cdot\bfU\bfW=\bfX^{-1}\bfQ_1^{(\alpha_1)}\bfY.
\end{align*}

\paragraph{Step 4 (circular permutation and remaining cores)}
Following \eqref{eq:Qk}, for any $\bgamma\in[n_1]\times\cdots\times[n_d]$,
\begin{align*}
    \hat{\bfQ}_{2\left<1\right>}&=\overleftarrow{\bfT}^{1}(:,\overleftarrow{\bgamma}_{\rm mid}^{1},\Gamma_1)\{\hat{\bfQ}_{1[1]}^*(:,\Gamma_1)\}^{-1}\\
    &=\left\{\bfQ_{2\left<1\right>}\left(\bfI_{r}\otimes(\bfQ_3^{(\gamma_3)}\cdots\bfQ_d^{(\gamma_d)})\right)\bfQ_{1[1]}^*(:,\Gamma_1)\right\}\left\{(\bfY^*\otimes\bfX^{-1})\bfQ_{1[1]}^*(:,\Gamma_1)\right\}^{-1}\\
    &=\bfQ_{2\left<1\right>}\left(\bfI_{r}\otimes(\bfQ_3^{(\gamma_3)}\cdots\bfQ_d^{(\gamma_d)})\right)\bfQ_{1[1]}^*(:,\Gamma_1)\{\bfQ_{1[1]}^*(:,\Gamma_1)\}^{-1}\left((\bfY^*)^{-1}\otimes\bfX\right)\\
    &=\bfQ_{2\langle1\rangle}\left((\bfY^*)^{-1}\otimes(\bfQ_3^{(\gamma_3)}\cdots\bfQ_d^{(\gamma_d)}\bfX)\right),
\end{align*}
which implies that $\hat{\bfQ}_2^{(\alpha_2)}=\bfY^{-1} \cdot \bfQ_2^{(\alpha_2)}\cdot (\bfQ_3^{(\gamma_3)}\cdots\bfQ_d^{(\gamma_d)}\bfX)$ for $\alpha_2\in[n_2]$. Similarly, for $k=3$,
\begin{align*}
    \hat{\bfQ}_{3\left<1\right>}&=\overleftarrow{\bfT}^2(:,\overleftarrow{\bgamma}_{\rm mid}^2,\Gamma_2)\left(\left(\bfI_{r}\otimes(\hat{\bfQ}_1^{(\gamma_1)})\right)\bfQ_{2[1]}^*(:,\Gamma_2)\right)^{-1}\\
    &=\left\{\bfQ_{3\left<1\right>}\left(\bfI_{r}\otimes(\bfQ_4^{(\gamma_4)}\cdots\bfQ_d^{(\gamma_d)}\bfQ_1^{(\gamma_1)})\right)\bfQ_{2[1]}^*(:,\Gamma_2)\right\} \\
    &\qquad\cdot\left\{\left(\bfI_{r}\otimes \hat{\bfQ}_1^{(\gamma_1)}\right)\left(((\bfQ_3^{(\gamma_3)}\cdots\bfQ_d^{(\gamma_d)}\bfX)^*\otimes\bfY^{-1}\right)\bfQ_{2[1]}^*(:,\Gamma_2)\right\}^{-1}\\
    &=\bfQ_{3\left<1\right>}\left(\bfI_{r}\otimes(\bfQ_4^{(\gamma_4)}\cdots\bfQ_d^{(\gamma_d)}\bfQ_1^{(\gamma_1)})\right)\left(((\bfQ_3^{(\gamma_3)}\cdots\bfQ_d^{(\gamma_d)}\bfX)^{-1})^*\otimes(\bfY(\hat{\bfQ}_1^{(\gamma_1)})^{-1})\right)\\
    &=\bfQ_{3\left<1\right>}\left(\bfI_{r}\otimes(\bfQ_4^{(\gamma_4)}\cdots\bfQ_d^{(\gamma_d)}\bfQ_1^{(\gamma_1)})\right)\left(((\bfQ_3^{(\gamma_3)}\cdots\bfQ_d^{(\gamma_d)}\bfX)^{-1})^*\otimes((\bfQ_1^{(\gamma_1)})^{-1}\bfX)\right)\\
    &=\bfQ_{3\left<1\right>}\left(((\bfQ_3^{(\gamma_3)}\cdots\bfQ_d^{(\gamma_d)}\bfX)^{-1})^*\otimes((\bfQ_4^{(\gamma_4)}\cdots\bfQ_d^{(\gamma_d)}\bfX)\right).
\end{align*}
Hence for $\alpha_3\in[n_3]$, $\hat{\bfQ}_3^{(\alpha_3)}=(\bfQ_3^{(\gamma_3)}\cdots\bfQ_d^{(\gamma_d)}\bfX)^{-1} \cdot\bfQ_3^{(\alpha_3)}\cdot(\bfQ_4^{(\gamma_4)}\cdots\bfQ_d^{(\gamma_d)}\bfX)$.
By induction, for $k\geq 3$, 
\begin{align*}
    \hat{\bfQ}_{k\left<1\right>}&=\overleftarrow{\bfT}^{k-1}(:,\overleftarrow{\bgamma}_{\rm mid}^{k-1},\Gamma_{k-1})\left(\left(\bfI_{r}\otimes(\hat{\bfQ}_1^{(\gamma_1)}\cdots\hat{\bfQ}_{k-2}^{(\gamma_{k-2})})\right)\hat{\bfQ}_{(k-1)[1]}^*(:,\Gamma_{k-1})\right)^{-1}\\
    &=\overleftarrow{\bfT}^{k-1}(:,\overleftarrow{\bgamma}_{\rm mid}^{k-1},\Gamma_{k-1})\left(\left(\bfI_{r}\otimes(\bfX^{-1}\bfQ_1^{(\gamma_1)}\cdots\bfQ_d^{(\gamma_d)}\bfX)\right)\hat{\bfQ}_{(k-1)[1]}^*(:,\Gamma_{k-1})\right)^{-1}\\
    &=\left\{\bfQ_{k\left<1\right>}\left(\bfI_{r}\otimes(\bfQ_{k+1}^{(\gamma_{k+1})}\cdots\bfQ_d^{(\gamma_d)}\bfQ_1^{(\gamma_1)}\cdots\bfQ_{k-2}^{(\gamma_{k-2})})\right)\right\}\\
    &\quad \cdot \left\{\left(\bfI_{r}\otimes(\bfX^{-1}\bfQ_1^{(\gamma_1)}\cdots\bfQ_d^{(\gamma_d)}\bfX)\right)\left((\bfQ_k^{(\gamma_k)}\cdots\bfQ_d^{(\gamma_d)}\bfX)^* \otimes (\bfQ_{k-1}^{(\gamma_{k-1})}\cdots\bfQ_d^{(\gamma_d)}\bfX)^{-1}\right)\right\}^{-1}\\
    &=\bfQ_{k\langle1\rangle}\left(((\bfQ_k^{(\gamma_k)}\cdots\bfQ_d^{(\gamma_d)}\bfX)^{-1})^* \otimes (\bfQ_{k+1}^{(\gamma_{k+1})}\cdots\bfQ_d^{(\gamma_d)}\bfX)\right).
\end{align*}
Therefore, 
\begin{align}\label{eq:q-equivalent}
    \hat{\mbcQ}_k=\left\{\begin{aligned}
        &\mbcQ_1 \times_2 \bfX^{-1} \times_3 \bfY^*, &&\text{if}\ k=1;\\
        &\mbcQ_2 \times_2 \bfY^{-1} \times_3 \{\bfQ_3^{(\gamma_3)}\cdots \bfQ_d^{(\gamma_d)}\bfX\}^*, &&\text{if}\ k=2;\\
        &\mbcQ_k \times_2 (\bfQ_k^{(\gamma_k)}\cdots\bfQ_d^{(\gamma_d)}\bfX)^{-1} \times_3 \{\bfQ_{k+1}^{(\gamma_{k+1})}\cdots\bfQ_d^{(\gamma_d)}\bfX\}^*, &&\text{if}\ 3\leq k\leq d-1;\\
        &\mbcQ_d \times_2 (\bfQ_d^{(\gamma_d)}\bfX)^{-1} \times_3 \bfX^*, &&\text{if}\ k=d;
    \end{aligned}\right.
\end{align}
which is exactly the equation~\eqref{eq:q-equivalent-main}.

\section{Proofs of Lemmas}\label{sec:proof-lemma}

\begin{proof}[Proof of Lemma~\ref{lemma:eigenvalues}]
    Let $\bfU=[\bfu_1\ \cdots\ \bfu_{r}]$. Then $\bfu_j\in\bbC^{r}$ is the eigenvector corresponding to $\lambda_j$ for $j\in[r]$. By \eqref{eq:trace} and the fact that with probability one, $\bfQ_{1\left<1\right>}$ is of full column rank, $\bfR^{\bbeta}\in\text{GL}(r,\bbC),$ and $\bfQ_{d[1]}^*(:,\Gamma_d)\in\text{GL}(r^2,\bbC)$, we have
    \begin{align*}
        \bfT(:,\balpha,\Gamma_{\balpha}) &= \bfQ_{1\left<1\right>}\left(\bfI_r\otimes \bfR^{\balpha}\right)\bfQ_{d[1]}^*(:,\Gamma_d),\quad \text{and}\\
        \bfT(:,\bbeta,\Gamma_{\bbeta})^\dagger&=\{\bfQ_{d[1]}^*(:,\Gamma_d)\}^{-1}\left(\bfI_{r}\otimes(\bfR^{\bbeta})^{-1}\right)\bfQ_{1\left<1\right>}^\dagger.
    \end{align*}
    Therefore,
    \begin{align*}
        \bfT(:,\balpha,\Gamma_{\balpha})\bfT(:,\bbeta,\Gamma_{\bbeta})^\dagger&=\bfQ_{1\left<1\right>}\left(\bfI_{r}\otimes\bfR^{\balpha}\right)\bfQ_{d[1]}^*(:,\Gamma_d)\{\bfQ_{d[1]}^*(:,\Gamma_d)\}^{-1}\left(\bfI_{r}\otimes(\bfR^{\bbeta})^{-1}\right)\bfQ_{1\left<1\right>}^\dagger\\
        &=\bfQ_{1\left<1\right>}\left(\bfI_{r}\otimes\left(\bfR^{\balpha}(\bfR^{\bbeta})^{-1}\right)\right)\bfQ_{1\left<1\right>}^\dagger\\
        &=\bfQ_{1\left<1\right>}\left(\bfI_{r}\otimes\left(\bfU\boldsymbol{\Lambda}\bfU^{-1}\right)\right)\bfQ_{1\left<1\right>}^\dagger\\
        &=\bfQ_{1\left<1\right>}\left(\bfI_{r}\otimes\bfU\right)\left(\bfI_{r}\otimes\boldsymbol{\Lambda}\right)\left(\bfI_{r}\otimes\bfU\right)^{-1}\bfQ_{1\left<1\right>}^\dagger.
    \end{align*}
    Hence the nonzero eigenvalues of $\bfT(:,\balpha,\Gamma_{\balpha})\bfT(:,\bbeta,\Gamma_{\bbeta})^\dagger$ are $\lambda_1,\dots,\lambda_{r}$, each with the geometric multiplicity $r$. For $k\in[r]$, the corresponding eigenspace is
    \begin{align*}
        \text{span}\left\{\bfQ_{1\left<1\right>}\begin{pmatrix}
            \bfu_k\\
            \mathbf{0}_{r(r-1)}
        \end{pmatrix},\dots,\bfQ_{1\left<1\right>}\begin{pmatrix}
            \mathbf{0}_{r(r-1)}\\
            \bfu_k
        \end{pmatrix}\right\}=\text{span}\left\{\bfQ_{1\left<1\right>}\left(\bfI_{r}\otimes\bfu_k\right)\bfK_k\right\},
    \end{align*}
    for all $\bfK_k\in\text{GL}(r,\bbC)$. Notice that for $\ell\in[n_1],\ j\in[r],\ k\in[r]$,
    \begin{align*}
        \bfE_{\ell,(j-1)r+k}&=\sum_{m=1}^{r^2}\left(\bfQ_{1\left<1\right>}\right)_{\ell,m}\left(\left(\bfI_{r}\otimes\bfU\right)\Pi(\bfK)\right)_{m,(j-1)r+k}\\
        &=\sum_{j_1=1}^{r}\sum_{k_1=1}^{r}\left(\bfQ_{1\left<1\right>}\right)_{\ell,(j_1-1)r+k_1} \left(\left(\bfI_{r}\otimes\bfU\right)\Pi(\bfK)\right)_{(j_1-1)r+k_1,(j-1)r+k}\\
        &=\sum_{j_1=1}^{r}\sum_{k_1=1}^{r}\left(\bfQ_{1\left<1\right>}\right)_{\ell,(j_1-1)r+k_1}\bfU_{k_1,k}\left(\bfK_k\right)_{j_1,j}\\
        &=\sum_{j_1=1}^{r}\sum_{k_1=1}^{r}\left(\bfQ_{1\left<1\right>}\right)_{\ell,(j_1-1)r+k_1}(\bfu_{k})_{k_1}\left(\bfK_k\right)_{j_1,j}\\
        &=\sum_{j_1=1}^{r}\sum_{k_1=1}^{r}\left(\bfQ_{1\left<1\right>}\right)_{\ell,(j_1-1)r+k_1}\left(\bfI_{r}\otimes\bfu_k\right)_{(j_1-1)r+k_1,j_1}\left(\bfK_k\right)_{j_1,j}\\
        &=\sum_{j_1=1}^{r}\sum_{k_1=1}^{r}\left(\bfQ_{1\left<1\right>}\right)_{\ell,(j_1-1)r+k_1}\left(\left(\bfI_{r}\otimes\bfu_k\right)\bfK_k\right)_{(j_1-1)r+k_1,j}\\
        &=\sum_{m=1}^{r^2}\left(\bfQ_{1\left<1\right>}\right)_{\ell,m}\left(\left(\bfI_{r}\otimes\bfu_k\right)\bfK_k\right)_{m,j}\\
        &=\left(\bfQ_{1\left<1\right>}\left(\bfI_{r}\otimes\bfu_k\right)\bfK_k\right)_{\ell,j}.
    \end{align*}
    Therefore, for $k\in[r]$, the space spanned by the set of the $\{(j-1)r+k\}$-column ($j\in[r]$) of $\bfE$ is exactly the eigenspace of $\bfT(:,\balpha,\Gamma_{\balpha})\bfT(:,\bbeta,\Gamma_{\bbeta})^{-1}$ corresponding to the $k$th eigenvalue $\lambda_k$.
\end{proof}

\begin{proof}[Proof of Lemma~\ref{lemma:hat-q-tilde}]
    By \eqref{eq:q-equivalent-main},
    \begin{align*}
        \hat{\bfQ}_k^{(\alpha_k)}=\bfX_k^{-1}\bfQ^{(\alpha_k)}\bfX_{k+1},\ k\in[d-1],\quad\text{and}\quad
        \hat{\bfQ}_d^{(\alpha_d)}=\bfX_d^{-1}\bfQ^{(\alpha_d)}\bfX_1,
    \end{align*}
    for some $\bfX_1,\bfX_2,\dots,\bfX_d\in\text{GL}(r,\bbC)$. Therefore,
    \begin{align*}
        \hat{\bfQ}_1^{(\alpha_1)}\hat{\bfQ}_2^{(\alpha_2)}\cdots\hat{\bfQ}_d^{(\alpha_d)}=\bfX_1^{-1}\bfQ^{(\alpha_1)}\bfQ^{(\alpha_2)}\cdots\bfQ^{(\alpha_d)}\bfX_1.
    \end{align*}
    Letting $\tilde{\bfQ}^{(\alpha_k)}=\bfX_1^{-1}\bfQ^{(\alpha_k)}\bfX_1$, $k\in[d]$ completes the proof. 
\end{proof}

\begin{proof}[Proof of Lemma~\ref{lemma:q}]
    By \eqref{eq:tilde.q}, we have
    \begin{align*}
        \bfY_j\boldsymbol{\Lambda}_j\bfY_j^{-1}=\hat{\bfQ}_1^{(j)}\hat{\bfQ}_2^{(j)}\cdots\hat{\bfQ}_d^{(j)}=(\tilde{\bfQ}^{(j)})^d.
    \end{align*}
    Therefore,
    \begin{align*}
        \boldsymbol{\Lambda}_j=\bfY_j^{-1}(\tilde{\bfQ}^{(j)})^d\bfY_j=\left(\bfY_j^{-1}\tilde{\bfQ}^{(j)}\bfY_j\right)^d,
    \end{align*}
    which implies that
    \begin{align*}
        \tilde{\bfQ}^{(j)}=\bfY_j\boldsymbol{\Omega}_j\bfY_j^{-1},
    \end{align*}
    for some $\boldsymbol{\Omega}_j=\text{diag}(\omega_{j,1},\dots,\omega_{j,r})$ such that $\omega_{j,k}^d=\lambda_{j,k}$, $k\in[r]$.
\end{proof}

\begin{proof}[Proof of Lemma~\ref{lemma:roots}]
    Assume that $\bar{\bfQ}^{(1)}=\bfY_1\boldsymbol{\Omega}_1\bfY_1^{-1}$ is a solution for \eqref{eq:tilde.q}, where $\boldsymbol{\Omega}_1=\text{diag}(\omega_{1,1},\dots,\omega_{1,r})$, and $\omega_{1,t}=\ell_t^{1/d}e^{i(\theta_t+2k_t\pi)/d}$ for $t\in[r]$. Then for $j\in\{2,\dots,d\}$,
    \begin{align*}
        \bar{\bfQ}^{(j)}=(\bar{\bfQ}^{(1)})^{-(d-1)}\hat{\bfQ}_1^{(1)}\hat{\bfQ}_d^{(j)}=\bfY_1\boldsymbol{\Omega}_1^{-(d-1)}\bfY_1^{-1}\hat{\bfQ}_1^{(1)}\hat{\bfQ}_d^{(j)}.
    \end{align*}
    Therefore, for all $\balpha\in[n]^d$, we have
    \begin{align}\label{eq:qj1-n}
        &\hat{\bfQ}_1^{(\alpha_1)}\hat{\bfQ}_2^{(\alpha_2)}\cdots\hat{\bfQ}_d^{(\alpha_d)}=\bar{\bfQ}^{(\alpha_1)}\bar{\bfQ}^{(\alpha_2)}\cdots\bar{\bfQ}^{(\alpha_d)}\nonumber\\
        &=\bfY_1\boldsymbol{\Omega}_1^{-(d-1)}\bfY_1^{-1}\hat{\bfQ}_1^{(1)}\hat{\bfQ}_d^{(\alpha_1)}\cdot\bfY_1\boldsymbol{\Omega}_1^{-(d-1)}\bfY_1^{-1}\hat{\bfQ}_1^{(1)}\hat{\bfQ}_d^{(\alpha_2)}\cdots\bfY_1\boldsymbol{\Omega}_1^{-(d-1)}\bfY_1^{-1}\hat{\bfQ}_1^{(1)}\hat{\bfQ}_d^{(\alpha_d)}.
    \end{align}
    Note that the above equation \eqref{eq:qj1-n} can be written element-wise as
    \begin{align}\label{eq:elementwise}
        \left(\hat{\bfQ}_1^{(\alpha_1)}\hat{\bfQ}_2^{(\alpha_2)}\cdots\hat{\bfQ}_d^{(\alpha_d)}\right)_{(j_1,j_2)}=\sum_{\btau\in\mcI(j_1,j_2)}c(j_1,j_2,\btau)\cdot\omega_{1,1}^{-(d-1)\tau_1}\omega_{1,2}^{-(d-1)\tau_2}\cdots\omega_{1,r}^{-(d-1)\tau_r},
    \end{align}
    where for $1\leq j_1,j_2\leq r$, $\mcI(j_1,j_2)$ is the set of some indices $\btau=(\tau_1,\tau_2,\dots,\tau_r)$ satisfying $0\leq\tau_1,\tau_2,\dots,\tau_r\leq d$ and $\tau_1+\tau_2+\cdots+\tau_r=d$, and $c(s_1,s_2,\btau)\in\bbC$ is the corresponding coefficient depending on the indices $(j_1,j_2)$ and $\btau\in\mcI(j_1,j_2)$. For $\omega_{1,t}'=\ell_t^{1/d}e^{(\theta_t+2(k_t+1)\pi)i/d}=\omega_{1,t}\cdot e^{2\pi i/d}$, $t\in[r]$, we have that
    \begin{align*}
        &\sum_{\btau\in\mcI(j_1,j_2)}c(j_1,j_2,\btau)\cdot(\omega_{1,1}')^{-(d-1)\tau_1}(\omega_{1,2}')^{-(d-1)\tau_2}\cdots(\omega_{1,r}')^{-(d-1)\tau_r}\\
        &=\sum_{\btau\in\mcI(j_1,j_2)}c(j_1,j_2,\btau)\cdot\left(\omega_{1,1}e^{2\pi i/d}\right)^{-(d-1)\tau_1}\left(\omega_{1,2}e^{2\pi i/d}\right)^{-(d-1)\tau_2}\cdots\left(\omega_{1,r}e^{2\pi i/d}\right)^{-(d-1)\tau_r}\\
        &=\sum_{\btau\in\mcI(j_1,j_2)}c(j_1,j_2,\btau)\cdot\omega_{1,1}^{-(d-1)\tau_1}\omega_{1,2}^{-(d-1)\tau_2}\cdots\omega_{1,r}^{-(d-1)\tau_r}\cdot(e^{2\pi i/d})^{-(d-1)\sum_{t=1}^r\tau_t}\\
        &=\sum_{\btau\in\mcI(j_1,j_2)}c(j_1,j_2,\btau)\cdot\omega_{1,1}^{-(d-1)\tau_1}\omega_{1,2}^{-(d-1)\tau_2}\cdots\omega_{1,r}^{-(d-1)\tau_r}\cdot(e^{2\pi i/d})^{-d(d-1)}\\
        &=\sum_{\btau\in\mcI(j_1,j_2)}c(j_1,j_2,\btau)\cdot\omega_{1,1}^{-(d-1)\tau_1}\omega_{1,2}^{-(d-1)\tau_2}\cdots\omega_{1,r}^{-(d-1)\tau_r}\\
        &=\left(\hat{\bfQ}_1^{(\alpha_1)}\hat{\bfQ}_2^{(\alpha_2)}\cdots\hat{\bfQ}_d^{(\alpha_d)}\right)_{(j_1,j_2)},
    \end{align*}
    where the last equality follows from \eqref{eq:elementwise}. Therefore, $\boldsymbol{\Omega}_1'=\text{diag}(\omega_{1,t}',\dots,\omega_{1,r}')$ is also a solution for \eqref{eq:tilde.q}.
\end{proof}

\begin{proof}[Proof of Lemma~\ref{lemma:mps}]
    By Lemma~\ref{lemma:permutation}, it suffices to show this for $j = d, k = 1$, in which case 
    $\tau_{\balpha'} = \bfT(1,\balpha',1)$ 
    and $\bfR^{d,1;\balpha} = \bfR^{\balpha}$.
    
    Contract $\brho^S$ in the mode corresponding to index $h$ along this unit vector to obtain the $n_1n_d\times n_1n_d$ matrix $\brho^{h,v}$ whose $((\alpha_1,\alpha_d), (\alpha'_1,\alpha'_d))$-th entry is given by
    \begin{equation}
        \brho^{h,v}(\alpha_1\alpha_d, \alpha'_1\alpha'_d) = \sum^{n_h}_{\alpha_h, \alpha'_h =1} v_{\alpha_h} v^*_{\alpha'_h} \brho^S(\alpha_1\alpha_h\alpha_d, \alpha'_1\alpha'_h\alpha'_d)\,.
    \end{equation}
    Analogously to Eq.~\eqref{eq:trace}, we can write the $n_1^2\times n_d^2$ reshaping of $\rho^{h,v}$ as
    \begin{equation}
        \label{eq:rho-reshaping}
        \sum_{\balpha_{\backslash h}} \sum^{n_h}_{\alpha_h, \alpha'_h = 1} v_{\alpha_h} v^*_{\alpha'_h} \left(\bfQ_{1\left<1\right>} (\bfI_r\otimes \bfR^{\balpha} ) \bfQ^*_{d[1]}\right) \otimes \left(\bfQ_{1\left<1\right>}  (\bfI_r\otimes \bfR^{\balpha'} ) \bfQ^*_{d[1]}\right)^{\mathrm{conj}}\,,
    \end{equation}
    where $\bfA^\mathrm{conj}$ denote the entrywise conjugate of a matrix $\bfA$.
    
    In particular, by looking at the $((\alpha_1,1),(\alpha_d,1))$-th entries of \eqref{eq:rho-reshaping} and scaling by a known real-valued constant factor, we obtain the submatrix $\bfM^{h,v}_{d,1}$ as claimed.
\end{proof}

\section{Refined Algorithm with Augmented Observations}\label{sec:algorithm}

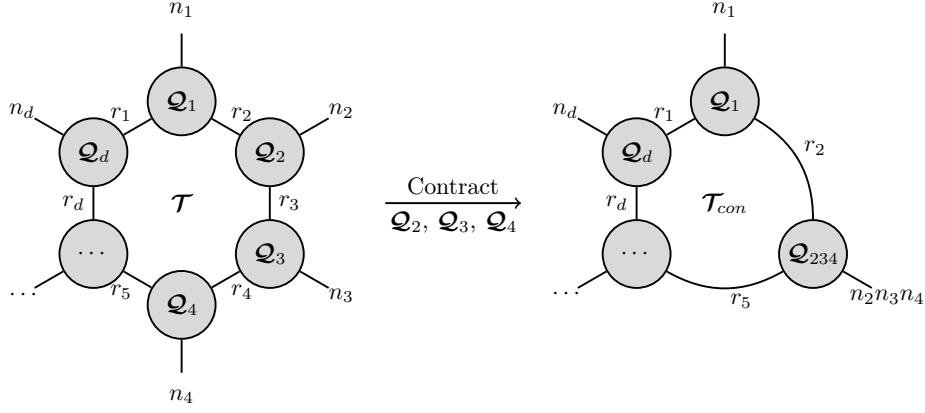
\begin{figure}[ht]
    \centering
    \begin{tikzpicture}[tensor/.style={circle, draw, fill=gray!30, minimum size=1cm, inner sep=0pt},
  connect/.style={-},
  thick]

\def\r{1.5}

\foreach \angle/\name/\label in {
  90/Z1/\mbcQ_1,
  30/Z2/\mbcQ_2,
 -30/Z3/\mbcQ_3,
 -90/Z4/\mbcQ_4
} {
  \node[tensor] (\name) at (\angle:\r) {$\label$};
}

\node[tensor] (Z5) at (-150:\r) {$\cdots$};
\node[tensor] (Z6) at (150:\r) {$\mbcQ_d$};
\node[] (T) at (0,0) {$\mbcT$};

\draw (Z1) -- node[midway, above right=-2pt] {$r_2$} (Z2);
\draw (Z2) -- node[midway, right] {$r_3$} (Z3);
\draw (Z3) -- node[midway, below right=-2pt] {$r_4$} (Z4);
\draw (Z4) -- node[midway, below left=-2pt] {$r_5$} (Z5);
\draw (Z5) -- node[midway, left] {$r_{d}$} (Z6);
\draw (Z6) -- node[midway, above left=-2pt] {$r_1$} (Z1);

\foreach \name/\angle in {Z1/90, Z2/30, Z3/-30, Z4/-90, Z5/-150, Z6/150} {
  \draw[connect] (\name) -- ++(\angle:1)
    node[pos=2.2, anchor=\angle] {
      $\ifnum\angle=90 n_1\else
       \ifnum\angle=30 n_2\else
       \ifnum\angle=-30 n_3\else
       \ifnum\angle=-90 n_4\else
       \ifnum\angle=-150 \cdots\else
       n_d\fi\fi\fi\fi\fi$
    };
}

\draw[->, thick] (3,0) to node[midway, above, black] {Contract} node [midway, below, black] {$\mbcQ_2$, $\mbcQ_3$, $\mbcQ_4$} (5,0);

\begin{scope}[xshift=8cm]

\node[tensor] (Z234) at (-30:\r) {$\mbcQ_{234}$};
\node[tensor] (Z5) at (-150:\r) {$\cdots$};
\node[tensor] (Z6) at (150:\r) {$\mbcQ_d$};
\node[tensor] (Z1) at (90:\r) {$\mbcQ_1$};
\node[] (Tcon) at (0,0) {$\mbcT_{con}$}; 

\draw[connect] (Z234) to [bend left=30] node [midway, below right=-2pt] {$r_5$} (Z5);
\draw[connect] (Z5) -- node [midway, left] {$r_d$} (Z6);
\draw[connect] (Z6) -- node [midway, above left=-2pt] {$r_1$} (Z1);
\draw[connect] (Z1) to [bend left=30] node [midway, above right=-2pt] {$r_2$} (Z234);

\draw[connect] (Z234) -- ++(-30:1) node[pos=2.5, anchor=-30] {$n_2n_3n_4$};
\draw[connect] (Z5) -- ++(-150:1) node[pos=2.2, anchor=-150] {$\cdots$};
\draw[connect] (Z6) -- ++(150:1) node[pos=2.2, anchor=150] {$n_d$};
\draw[connect] (Z1) -- ++(90:1) node[pos=2.2, anchor=90] {$n_1$};

\end{scope}

\end{tikzpicture}
    \caption{Contraction of TR cores.}
    \label{fig:contraction}
\end{figure}

In this section, we present a refined algorithm for exact TR decomposition \eqref{eq:order4TR} under a relaxed version of the conditions in Theorem~\ref{thm:order-d}. Specifically, we assume that all conditions hold except for $n_j \geq r^2$ and $n_{j+1} \geq r^2$ for some $j \in [d-1]$, and $n_k < r^2$ for some $k \in [d]$, while assuming access to additional observations. Let $A\subseteq[d]$ be the set of indices satisfying the dimension condition: 
\begin{align*}
    A=\{k\in[d]:n_k\geq r^2\}. 
\end{align*}
For $k\in A$, we define the function 
\begin{align*}
    \lambda(k)=\sup\left\{\ell \in [d]:\ (k+j-1)\bmod d \in A\ \text{for all}\ j\in[\ell]\right\},
\end{align*}
which represents the maximal length of a consecutive sequence (modulo $d$) contained in $A$, starting from the index $k$. Let $q=\argmax_{k\in A}\lambda(k)$ be the starting index of the longest such subsequence. The refined algorithm then begins by recovering the cores with indices $[(q,q+1,\dots,q+\lambda(q)-1)\bmod d]$ through \emph{contraction} of the remaining cores (Figure~\ref{fig:contraction}), and then expands the contracted representation to reconstruct the full set of TR cores $\{\hat{\mbcQ}_k\}_{k\in[d]}$.

Without loss of generality, we assume in what follows that $d=\argmax_{k\in A}\lambda(k)$. We then consider the following contracted TR decomposition
\begin{align}\label{eq:contractTR}
    T_{con}(\alpha_1,\dots,\alpha_{\lambda(d)-1},\overline{\balpha},\alpha_d)=\text{tr}\left\{\bfQ_1^{(\alpha_1)}\cdots\bfQ_{\lambda(d)-1}^{(\alpha_{\lambda(d)-1})}\bfQ_{con}^{\balpha}\bfQ_d^{(\alpha_d)}\right\},
\end{align}
where $\balpha=(\alpha_{\lambda(d)},\dots,\alpha_{d-1})\in[n_{\lambda(d)}]\times\cdots\times[n_{d-1}]$ for $\lambda(d)\leq d-1$ (In Algorithm~\ref{alg:order-d}, one considers the case $\lambda(d)=\argmax_{k\in A}\lambda(k)=d$); $\overline{\balpha}$ is the column-major vectorization of $\balpha$; $\mbcT_{con}$ is the contraction of tensor $\mbcT$ satisfying
\begin{align}\label{eq:contract}
    T_{con}(\alpha_1,\dots,\alpha_{\lambda(d)-1},\overline{\balpha},\alpha_d)=T(\alpha_1,\dots,\alpha_{\lambda(d)-1},\balpha,\alpha_d);
\end{align}
and $\bfQ_{con}^{\balpha}=\bfQ_{\lambda(d)}^{(\alpha_{\lambda(d)})}\cdots\bfQ_{d-1}^{(\alpha_{d-1})}$. Since the contracted representation \eqref{eq:contractTR} satisfies all the conditions in Theorem~\ref{thm:order-d}, one can implement Algorithm~\ref{alg:order-d} to obtain a set of representatives $\{\hat{\mbcQ}_1,\dots,\hat{\mbcQ}_{\lambda(d)-1},\hat{\mbcQ}_{con},\hat{\mbcQ}_d\}\sim\{\mbcQ_1,\dots,\mbcQ_{\lambda(d)-1},\mbcQ_{con},\hat{\mbcQ}_d\}$. Then decompose $\hat{\mbcQ}_{con}$ to get $\hat{\mbcQ}_{\lambda(d)},\dots,\hat{\mbcQ}_{d-1}$. Following a similar argument as in the proof of Theorem~\ref{thm:order-d} (Section~\ref{sec:proof}), one can show that the cores \(\{\hat{\mbcQ}_k\}_{k \in [d]}\) obtained from Algorithm~\ref{alg:refined} also satisfy \eqref{eq:q-equivalent}, and therefore constitute valid TR cores.

We summarize the procedure in Algorithm~\ref{alg:refined}. The entries required for recovering the TR decomposition \eqref{eq:order4TR} are 
\begin{equation}
\begin{split}
    \label{eq:Delta-rf}
    \Delta_{rf}\equiv  ~~ & \bigl\{\bfT_{con}(:,\balpha,\Gamma_{\balpha}),\,\bfT_{con}(:,\bbeta,\Gamma_{\bbeta})\bigr\}\, \bigcup\, \bigl\{\bfT_{con}(:,\balpha',\Gamma_{\balpha'}),\,\bfT_{con}(:,\bbeta',\Gamma_{\bbeta'})\bigr\}\\
    & 
    \,\bigcup\, \bigl\{\overleftarrow{\bfT_{con}}^k(:,\overleftarrow{\bgamma}^k,\Gamma_{\lambda(d)-1})\bigr\}_{k\in[\lambda(d)+1]}.
\end{split}
\end{equation}

\begin{algorithm}[ht]
    \caption{Order-$d$ TR decomposition with augmented observations}
    \label{alg:refined}
    \KwIn{Tensor $\mbcT\in\bbC^{n_1\times n_2\times\cdots\times n_d}$ observed at entries $\Delta_{rf}$ in \eqref{eq:Delta-rf}, TR-rank $r$, $\bgamma\in[n_1]\times\cdots\times[n_{\lambda(d)-1}]\times[n_{\lambda(d)}\cdots n_{d-1}]\times [n_d]$}
    \KwOut{TR-cores $\hat{\mbcQ}_k\in\bbC^{n_k\times r^2}$, $k\in[d]$}

    Contract $\mbcT$ to obtain $\mbcT_{con}$ by \eqref{eq:contract}\;
    
    Follow lines 1--6 in Algorithm~\ref{alg:order-d} with $\mbcT_{con}$ to obtain $\hat{\mbcQ}_1$\;
    
    \For {$2\leq k\leq \lambda(d)+1$}{ 
        Compute 
        \begin{align*}
            \hat{\bfQ}_{k\left<1\right>}=
                \overleftarrow{\bfT_{con}}^{k-1}(:,\overleftarrow{\bgamma}^{k-1},\Gamma_{\lambda(d)-1})\left(\left(\bfI_{r_2}\otimes(\hat{\bfQ}_1^{(1)}\cdots\hat{\bfQ}_{k-2}^{(1)})\right)\hat{\bfQ}_{(k-1)[1]}^*(:,\Gamma_{k-1})\right)^\dagger,
        \end{align*}
        with the convention that the product in the parentheses is $\bfI_r$ when $k=2$\;
    }\tcp{Where it is understood that  $\hat{\mbcQ}_{\lambda(d)}=\hat{\mbcQ}_{con}$ in this loop}
    \For {$\lambda(d)\leq k\leq d-1$}{
        \For {$\alpha_k\in[n_k]$}{
            Let $\hat{\bfQ}_k^{(\alpha_k)}=\hat{\bfQ}_{con}^{\mathbf{1}_{k,\alpha_k}}$, where 
            $\mathbf{1}_{k,\alpha_k}=(\underbrace{1,\dots,1}_{k-\lambda(d)},\ \alpha_k,\ \underbrace{1,\dots,1}_{d-k-1})$
        }
    }
\end{algorithm}

\section{Additional Numerical Results}\label{sec:addsimulation}

This section reports additional numerical results for Algorithm~\ref{alg:noise}. We evaluate the relative reconstruction error, $\|\mbcT_n - \hat{\mbcT}\|_F / \|\mbcT_n\|_F$ over a given number of iterations among 100 independent trials. Here, $\mbcT_n$ denotes a tensor with TR rank $r \in \{2,3,4,5\}$. In all settings, the entries of the TR cores $\{\mbcQ_k\}_{k\in[d]}$ are drawn independently from $\mcN(0,\sigma_s^2)$, and Gaussian noise is added with variance $\sigma_n^2$. The average relative errors after the fixed number of iterations are plotted against the noise scale $\sigma_n$.

\begin{figure}[ht]
    \centering
    \includegraphics[width=\linewidth]{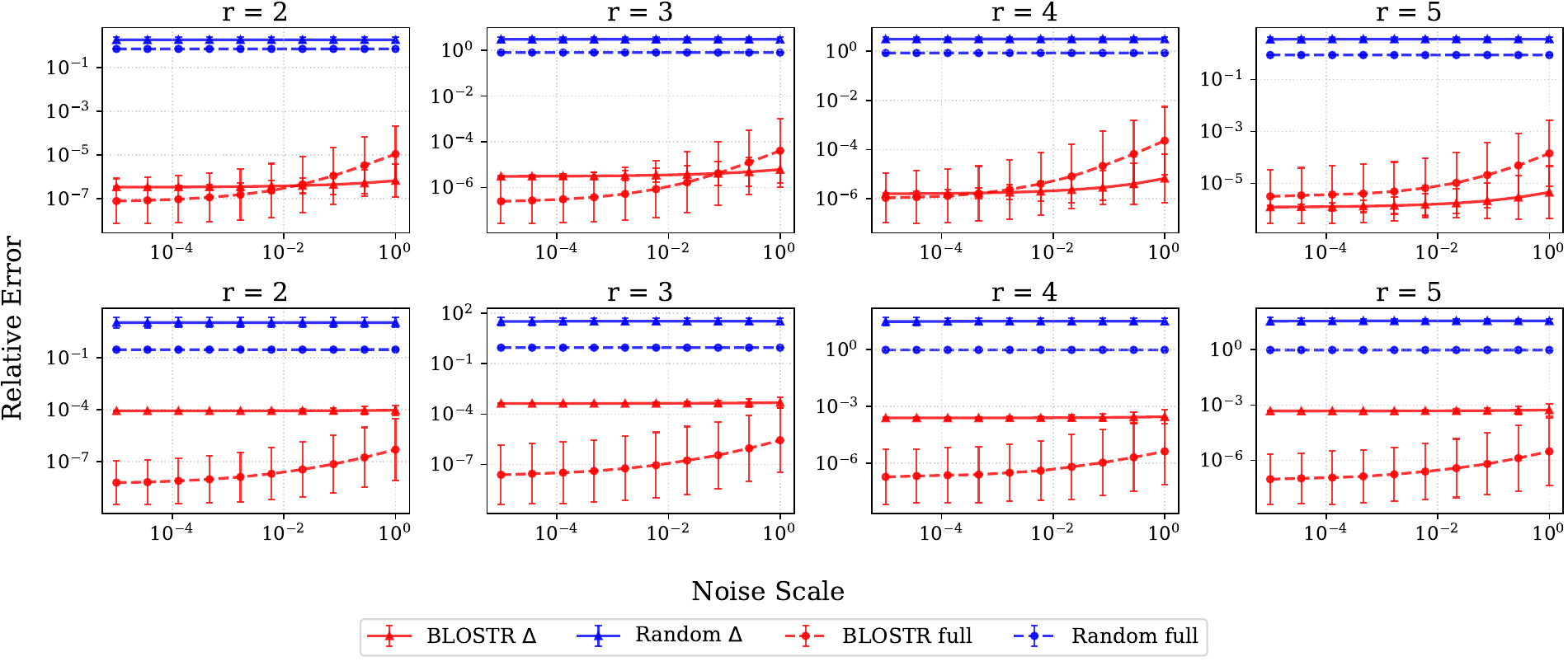}
    \caption{Average relative error (with log-scale standard deviation) of \blostr (Algorithm~\ref{alg:noise}) and randomly initialized ALS over increasing noise within 1 iteration. All other settings are the same as in Figure~\ref{fig:err-vs-ns-t3}.}
    \label{fig:err-vs-ns-t1}
\end{figure}

\begin{figure}[ht]
    \centering
    \includegraphics[width=\linewidth]{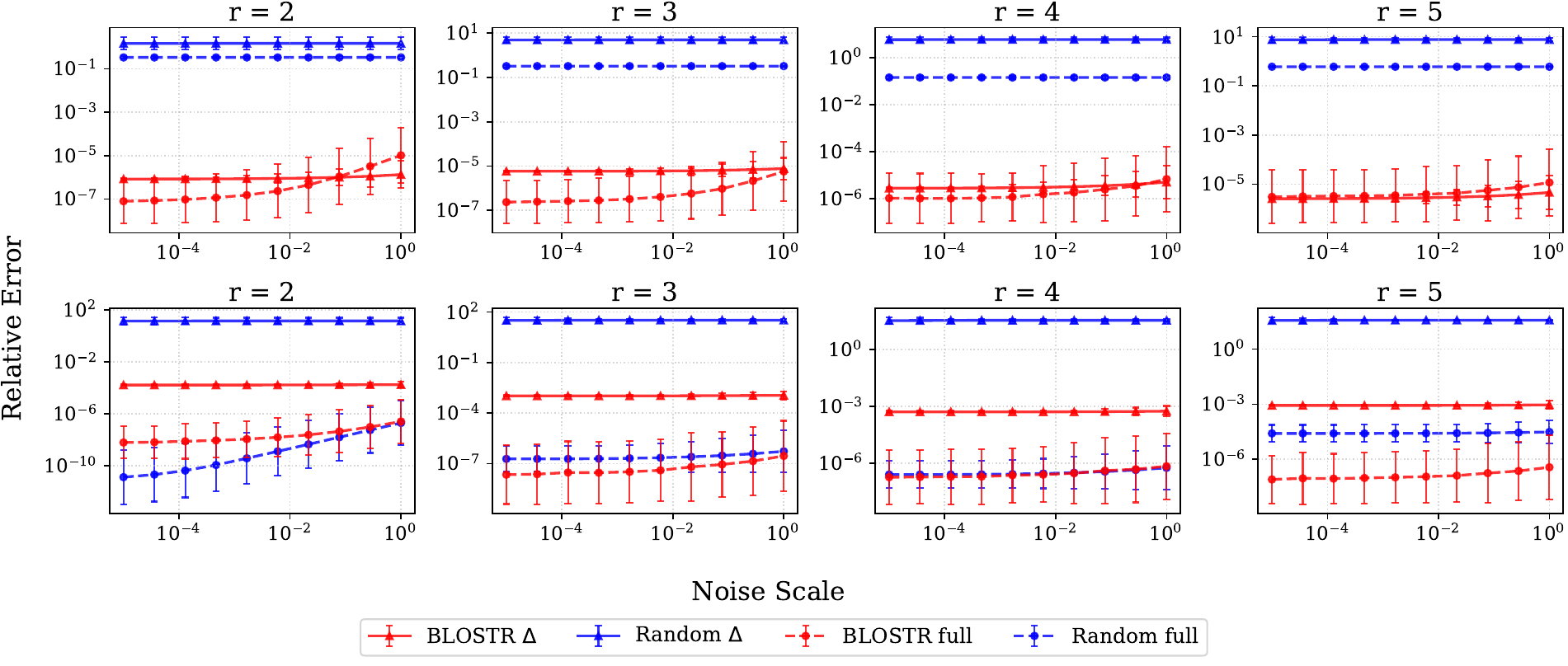}
    \caption{Average relative error (with log-scale standard deviation) of \blostr (Algorithm~\ref{alg:noise}) and randomly initialized ALS over increasing noise within 7 iterations. All other settings are the same as in Figure~\ref{fig:err-vs-ns-t3}.}
    \label{fig:err-vs-ns-t7}
\end{figure}

\begin{figure}[ht]
    \centering
    \includegraphics[width=\linewidth]{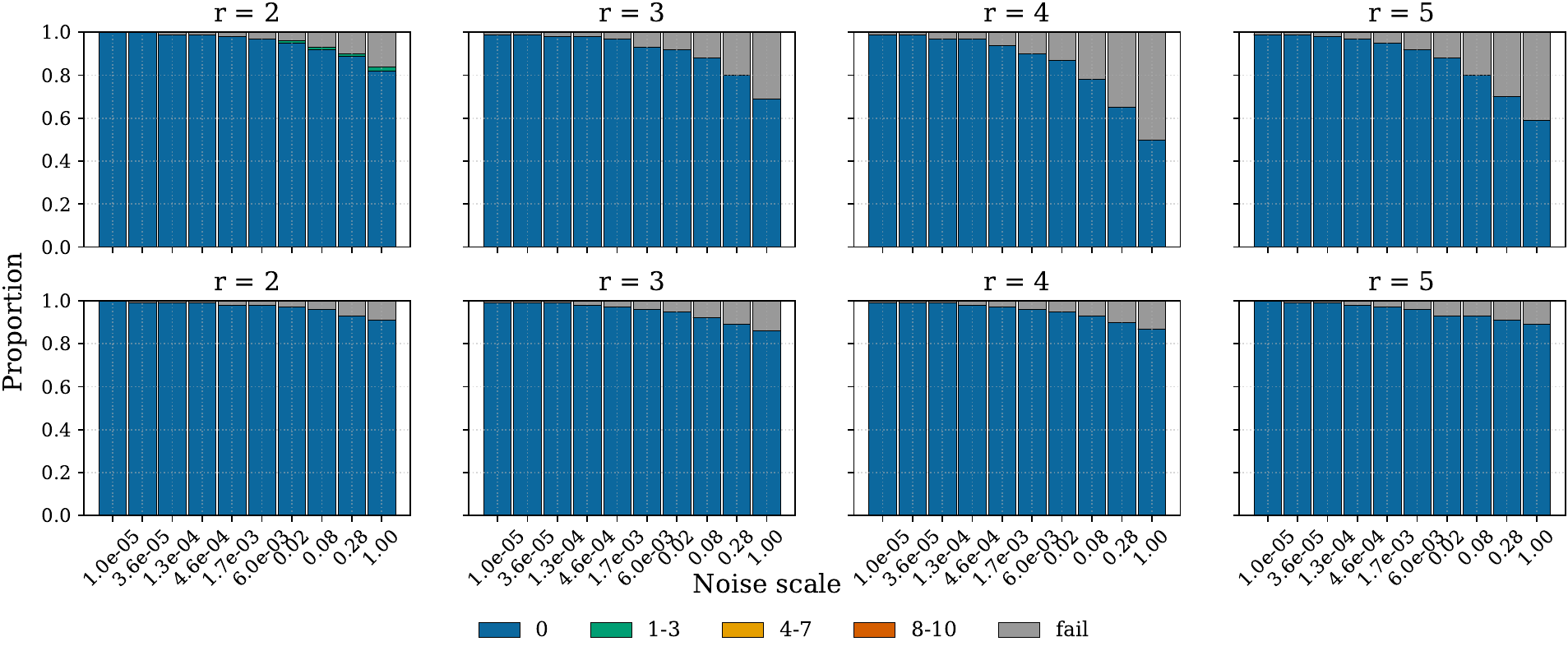}
    \caption{Proportion of successful recoveries (relative error below $10^{-6}$) within 10 iterations, averaged over 100 independent trials, using \blostr (Algorithm \ref{alg:noise}) under varying noise scales $\sigma_n$. Each tensor core $\{\mbcQ_k\}_{k\in[d]}$ has entries independently drawn from $\mcN(0, 10^2)$. The underlying tensor $\mbcT$ has dimension $30^{\times 3}$ in the top row and $30^{\times 4}$ in the bottom row.}
    \label{fig:t2tol-1e-6}
\end{figure}

\end{document}